\documentclass[12pt]{llncs}
\usepackage[utf8]{inputenc}
\usepackage[T1]{fontenc}
\usepackage{hyperref}
\usepackage{url}
\usepackage{booktabs}
\usepackage{amsfonts}
\usepackage{nicefrac}
\usepackage{microtype}
\usepackage{lipsum}
\usepackage{amsmath}
\usepackage{csvsimple}
\usepackage{tikz}
\usepackage{pgfplots}
\usepackage{pgfplotstable}
\usepackage{algorithm}
\usepackage{xcolor}
\usepackage[noend]{algpseudocode}
\usepackage{enumerate}
\usepackage{todonotes}
\usepackage[round,sectionbib]{natbib}
\usepackage[margin=1in]{geometry}

\allowdisplaybreaks
\newcommand{\TABLE}[3]{\centering\caption{#1}#2}
\newcommand{\FIGURE}[3]{\centering \protect{#1} \caption{#2}}
\renewcommand{\citep}[1]{(\cite{#1})}

\title{Learning to Solve Large-Scale Security-Constrained Unit Commitment Problems}

\author{
  \'Alinson S. Xavier \inst{1} \and
  Feng Qiu \inst{1} \and
  Shabbir Ahmed \inst{2}}

\institute{
  Energy Systems Division,
  Argonne National Laboratory,
  Argonne, IL, USA.  \email{\{axavier,fqiu\}@anl.gov}
  \and
  School of Industrial and Systems Engineering,
  Georgia Institute of Technology,
  Atlanta, GA, USA. \email{sahmed@isye.gatech.edu}
}

\begin{document}
\maketitle

\begin{abstract}
Security-Constrained Unit Commitment (SCUC) is a fundamental problem in power systems and electricity markets. In practical settings, SCUC is repeatedly solved via Mixed-Integer Linear Programming, sometimes multiple times per day, with only minor changes in input data. In this work, we propose a number of machine learning (ML) techniques to effectively extract information from previously solved instances in order to significantly improve the computational performance of MIP solvers when solving similar instances in the future. Based on statistical data, we predict redundant constraints in the formulation, good initial feasible solutions and affine subspaces where the optimal solution is likely to lie, leading to significant reduction in problem size. Computational results on a diverse set of realistic and large-scale instances show that, using the proposed techniques, SCUC can be solved on average 4.3x faster with optimality guarantees, and 10.2x faster without optimality guarantees, but with no observed reduction in solution quality. Out-of-distribution experiments provides evidence that the method is somewhat robust against dataset shift.
\end{abstract}

\section{Introduction}

Security-Constrained Unit Commitment (SCUC) is one of the most fundamental optimization problems in power systems and electricity markets, being solved a number of times daily by major reliability coordinators, independent system operators (ISOs), and regional transmission organizations (RTOs) to schedule 3.8 trillion kWh of energy production and clear a \$400 billion electricity market annually in the United States \citep{EiaWebsite}.
The problem asks for the most cost-efficient schedule for power generation units under a large number of physical, operational and economic constraints. Among other requirements, a feasible generation schedule must not only guarantee that enough power is being produced during each hour to satisfy the demand, but also that no transmission lines are overloaded during the delivery of the electric power. Furthermore, the schedule must be robust against small variations in demands and against the unexpected loss of a small number of network elements, such as generation units and transmission lines.

The computational performance of SCUC is an extremely important practical issue, given the very short clearing windows enforced in most electricity markets. For instance, in Midcontinent Independent System Operator (or MISO, one of the largest electricity markets in North America), systems operators only have 3 or 4 hours, after receiving all bids and offers, to produce a low-cost and fair generation schedule. During this short time window, SCUC must be solved multiple times, under a number of different scenarios. Using state-of-art software, realistic large-scale instances can be usually solved to 0.1\% optimality within approximately 20 minutes of running time \citep{ChenCastoWangWangWangWan2016}.
Improvements in computational performance would allow markets to implement numerous enhancements that could bring significant economic benefits and improved market efficiency, such as more accurate modelling of enhanced combined cycle units, higher resolution for production cost curves, sub-hourly dispatch, longer planning horizons, among others. Alternatively, a simple reduction of the optimality gap to 0.05\% without sacrificing computational running times could lead to significant cost savings, considering the \$400 billion market size in U.S. alone.

In the past, SCUC has been solved using various optimization techniques including priority lists \citep{BurnsGibson1975}, dynamic programming \citep{Lowery1966} and Lagrangian relaxation \citep{MerlinSandrin1983}. Nowadays, SCUC is most commonly formulated as a Mixed-Integer Linear Programming (MIP).
The introduction of MIP solvers allowed system operators to
easily model
more realistic constraints that had been traditionally hard to enforce with previous methods, and have benefited from the constantly advancing capabilities of modern solvers.
Since the first MIP formulations, proposed in the 1960s \citep{Garver1962}, a number of alternative formulations and polyhedral results have been described \citep{AtakanLulliSen2018,OstrowskiAnjosVannelli2012,RajanTakriti2005,GentileMoralesEspanaRamos2017,DamciKurtKucukyavuzRajanAtamturk2016,espana2015tight, lee2004min,pan2016polyhedral,Knueven2018}. There has also been active research on more practical techniques for improving the computational performance of the MIP approach, such as fast identification of redundant constraints \citep{ZhaiGuanChengWu2010,ArdakaniBouffard2015,XavierFengWangThimmapuram2018}, decomposition methods \citep{FeizollahiCostleyAhmedGrijalva2015, Kim18TemporalDecomp}, among others.

One aspect of SCUC that has been overlooked in previous research is that, in most practical settings, this problem is repeatedly solved with only minor changes in input data, often multiple times per day. Most of the system characteristics, such as the parameters describing the generation units or the topology of the transmission network, are kept almost entirely unchanged from solve to solve.
While modern MIP solvers have the ability to reuse previous solutions as warm starts, it is well known that, for SCUC, providing the previous-day optimal solution brings only minor performance improvements \citep{ChenCastoWangWangWangWan2016}.
In this work, we propose the use of machine learning to effectively extract information from previously solved instances, and show how to leverage this information to significantly accelerate the solution of similar instances in the future.

Machine learning (ML) is a discipline that studies how algorithms can automatically learn from large amounts of data and subsequently use the acquired knowledge to make decisions quickly and accurately.
In the 1990s and early 2000s, researchers experimented using artificial neural networks (ANN) to solve simpler variants of SCUC on small-scale instances \citep{SasakiWatanabeKubokawaYorinoYokoyama1992,WangShahidehpour1993,WalshOMalley1997,LiangKang2000}. Even when considering very simplified versions of the problem, obtaining sufficiently high-quality solutions to SCUC via ANN proved very challenging, and the approach failed to replace existing methods.
Similar explorations with evolutionary algorithms proved equally challenging \citep{JusteKitaTanakaHasegawa1999}.

In more recent years, there has been a growing interest, within the broader mathematical optimization community, in applying ML techniques to enhance the performance of current MIP solvers.
ML has been used, for example, to automatically construct branching strategies \citep{AlvarezLouveauxWehenkel2017,KhalilLeBodicSongNemhauserDilkina2016, LodiGiulia2017},
to better parallelize branch-and-bound methods \citep{AlvarezWehenkelLouveaux2014}, to decide when to run primal heuristics \citep{KhalilDilkinaNemhauserAhmedShao2017}, to predict resolution outcome \citep{FischettiLodiZarpellon2019}, to decide solver parameters \citep{BonamiLodiZarpellon2018} and
to automatically construct decompositions \citep{BassoCeselliTettamanzi2017}.
More broadly, there has been a renewed interest in using ML to tackle hard combinatorial problems. To give a few examples, deep and reinforcement learning have been used to produce heuristic solutions for optimization problems over graphs \citep{KhalilDaiZhangDilkina2017, kool2018attention, MichelPierreAlexandreYossiriLouisMartin2018} and for stochastic integer programs \citep{nair2018learning}. We refer to \citep{bengio2018machine} for a more complete survey. From a more abstract perspective, there has been research on integrating the fields of predictive, prescriptive and descriptive analytics \citep{lombardi2017empirical, bertsimas2019predictive}.
More directly related to our work, statistical learning has also been applied to the DC Optimal Power Flow Problem \citep{misra2018learning}.

In this work, instead of completely replacing existing MIP methods by ML models, as done in previous research with artificial neural networks and genetic algorithms, we propose the usage of ML techniques to significantly enhance the warm-start capabilities of modern MIP solvers when solving SCUC.
More specifically, we develop three ML models that can automatically identify different types of patterns in previously solved instances and subsequentially use these patterns to improve MIP performance.
The first model is designed to predict, based on statistical data, which constraints are necessary in the formulation and which constraints can be safely omitted. The second model constructs, based on a large number of previously obtained optimal solutions, a (partial) solution which is likely to work well as warm start. The third model identifies, with very high confidence, a smaller-dimensional affine subspace where the optimal solution is likely to lie, leading to the elimination of a large number of decision variables and significantly reducing the complexity of the problem.

We start, in Section \ref{sec:background}, by presenting the MIP formulation of SCUC and some brief introduction to the ML concepts we use.
In Section \ref{sec:learning}, we formalize the setting in which the learning takes place and we describe the proposed learning strategies. In Section \ref{sec:comp}, we evaluate the practical impact of the proposed ML approach by performing extensive computational experiments on a diverse set of large-scale and realistic SCUC instances, containing up to 6515 buses and 1388 units, under realistic uncertainty scenarios. Using the proposed technique, we show that SCUC can be solved on average 4.3x faster than conventional methods with optimality guarantees and 10.2x faster without optimality guarantees, but no observed reduction in solution quality.
We also conduct some experiments where the training and test samples are drawn from similar, but not exactly the same distribution, to provide some evidence that the method is somewhat robust against moderate dataset shifts.
Finally, in Section \ref{sec:limitations}, we discuss some limitations of the method, as well as directions for future work.
Although the techniques are presented in the context of one particular application, they can be easily adapted to other challenging optimization problems.

\section{Background} \label{sec:background}

~

\subsection{Mixed-Integer Linear Programming Formulation of SCUC}
\label{subsec:scuc}

Unit Commitment (UC) refers to a broad class of optimization problems dealing with the scheduling of power generating units \citep{CohenSherkat87}. For each hour in the planning horizon, the problem is to decide which units should be operational and how much power should they produce \citep{OstrowskiAnjosVannelli2012}. Each unit typically has minimum and maximum power outputs, a quadratic production cost curve, and fixed startup and shutdown costs. Other constraints typically include ramping rates, which restrict the maximum variation in power production from hour to hour, and minimum-up and minimum-down constraints, which prevents units from starting and shutting down too frequently.
The Security-Constrained Unit Commitment (SCUC) problem is a subclass of UC that, in addition to all the constraints previously mentioned, also guarantees the deliverability of power, by enforcing that the power flow in each transmission line does not exceed its safe operational limits \citep{Shaw95}.
To guarantee sufficient security margin in case of component failure, the problem considers not only transmission constraints under normal operating conditions, but also when there is one transmission line failure in the system (so-called N-1 transmission contingency) \citep{BatutRenaud92}. In this way, even if a transmission line unexpectedly fails and the power flow in the network changes, there will be no violations of the transmission line thermal limits, which is a reliability requirement enforced by North American Electric Reliability Corporation \citep{NERCwebsite}.

A number of MIP formulations and strong valid inequalities for SCUC have been proposed. In this work, we use the formulation proposed in \citep{MoralesEspanaLatorreRamos2013}, since it presents good computational performance even when transmission and N-1 security constraints are present. The full MIP formulation is shown in Appendix A. Here, we present a simplified summary, containing key decision variables and constraints that are most influential to the computational performance of the problem. We note, however, that the methods we present later sections make very few assumptions about the formulation, and can still be used with additional decision variables and constraints.

Consider a power system composed by a set $B$ of buses, a set $G$ of generators and a set $L$ of transmission lines. Let $T$ be the set of hours within the planning horizon. For each generator $g \in G$ and hour $t \in T$, let $x_{gt}$ be a binary decision variable indicating if the generator is operational, and $y_{gt}$ be a continuous decision variable indicating the amount of power being produced. The problem can then be formulated as
\begin{align}
  \text{minimize}
  \hspace{1em}& \sum_{g \in G} c_g(x_{g\bullet}, y_{g\bullet})
  \label{c:obj}\\
\text{subject to} \hspace{1em}
    & (x_{g\bullet}, y_{g\bullet}) \in \mathcal{G}_g & \forall g \in G
    \label{c:generator} \\
    &\sum_{g \in G} y_{gt} = \sum_{b \in B} d_{bt} & \forall b \in B, t \in T
    \label{c:balance} \\
    & -F^0_l \leq \sum_{b \in B} \delta^0_{lb} \left( \sum_{g \in G_b} y_{gt} - d_{bt} \right) \leq F^0_l & \forall l \in L, t \in T.
    \label{c:transmission} \\
    & -F^c_l \leq \sum_{b \in B} \delta^c_{lb} \left( \sum_{g \in G_b} y_{gt} - d_{bt} \right) \leq F^c_l & \forall c \in L, l \in L, t \in T.
    \label{c:security} \\
    & x_{gt} \in \{0,1\} & \forall g \in G, t \in T \\
    & y_{gt} \geq 0 & \forall g \in G, t \in T
\end{align}
In the objective function, $c_g$ is a piecewise-linear function which includes the startup, shutdown and production costs for generator $g$. The notation $(x_{g\bullet}, y_{g\bullet})$ is a shorthand for the vector $(x_{gt_1},\ldots,x_{gt_k}, y_{gt_1}, \ldots, y_{gt_k})$, where $\{t_1,\ldots,t_k\} = T$.
Constraint \eqref{c:generator} enforces a number of constraints which are local for each generator, including production limits, minimum-up and minimum-down running times, ramping rates, among others.
Constraint \eqref{c:balance} enforces that the total amount of power produced equals the sum of the demands $d_{bt}$ at each bus.
Constraints \eqref{c:transmission} enforces the deliverability of power under normal system conditions. The set $G_b$ denotes the set of generators attached to bus $b$.
We recall that power flows are governed by physical laws, i.e., Ohm's laws and Kirchhoff's Laws. Using the DC linearization of these laws, it is possible to express the thermal limits for each transmission line as \eqref{c:transmission}, where $\delta^0_{bl}$ are real numbers known as Injection Shift Factors (ISFs).
Similarly, constraints \eqref{c:security} enforce the deliverability of power under the scenario that line $c$ has suffered an outage. Note that the bounds and the ISFs may be different under each outage scenario.

To improve the strength of the linear relaxation, the full formulation contains additional auxiliary binary decision variables. However, once $x_{gt}$ is determined, all the other binary variables are immediately implied. Furthermore, once all the $x_{gt}$ variables are fixed, the problem reduces to a linear programming problem (known as Economic Dispatch) which can be solved efficiently.

Constraints \eqref{c:transmission} and \eqref{c:security} have a significant impact on the computational performance of SCUC.
The total number of such constraints is quadratic on the number of transmission lines. Since large power systems can have more than 10,000 lines, the total number of constraints can easily exceed hundreds of millions. To make matters worse, these constraints are typically very dense, which makes adding them all to the formulation impractical for large-scale instances, not only to due to degraded MIP performance, but even due to memory requirements.
Fortunately, it has been observed that enforcing only a very small subset of these constraints is already sufficient to guarantee that all the remaining ones are automatically satisfied \citep{BouffardGalianaArroyo2005}.
Identifying this critical subset of constraints can be very challenging, and usually involves either solving a relaxation of SCUC multiple times \citep{TejadaSanchezRamos2018,XavierFengWangThimmapuram2018} or solving auxiliary optimization problems \citep{ArdakaniBouffard2015,ZhaiGuanChengWu2010}.
In practice, system operators still rely on experience and external processes to pre-select a subset of constraints to include in the optimization problem.


\subsection{Machine Learning Classifiers}

In this work, we employ two classical supervised learning methods for binary classification: \emph{$k$-Nearest Neighbors} ($k$NN) and \emph{Support Vector Machines} (SVM). Given a training dataset containing labeled examples $(p^1,h_1),\ldots,(p^s,h_s) \in \mathbb{R}^d \times \{-1,1\}$, each method constructs a function $\phi:\mathbb{R}^d \to \{-1,1\}$ which can be used to label new samples $\tilde{p} \in \mathbb{R}^d$. In the following, we present a brief description of each method. For a more complete description, we refer to \cite{Alpaydin2014}.

$k$NN works by identifying the most similar examples in the training dataset and conducting a simple majority vote. More precisely, given a point $\tilde{p} \in \mathbb{R}^d$, the method first sorts the labeled examples by increasing distance $\|p^i - \tilde{p}\|$, according to some norm, to obtain a sequence $(p^{i_1},h_{i_1}),\ldots,(p^{i_s},h_{i_s})$. Then, it assigns label $1$ to $\tilde{p}$ if $\frac{h_{i_1} + \ldots h_{i_k}}{k} > 0$, and label $-1$ otherwise. More efficient implementations, which avoid sorting the training dataset for each evaluation, have also been described.
In Subsections~\ref{subsec:tr} and \ref{subsec:ws}, we use variations of this method for prediction of redundant constraints and warm starts.

Support Vector Machines construct a hyperplane that best separates training examples into two convex regions, according to their labels, and use this hyperplane for further classification. More precisely, the method first finds a hyperplane $(\pi,\pi_0) \in \mathbb{R}^d \times \mathbb{R}$ by solving the quadratic optimization problem
\begin{align*}
    \underset{\pi, \pi_0, \alpha}{\text{minimize }} \;\;
        & \frac{1}{2} \pi^T \pi + \frac{C}{s} \sum_{i=1}^s \alpha_i \\
    \text{subject to } \;\;
        & h_i \left( \pi^T p^i - \pi_0 \right) \geq 1 - \alpha_i
            & i=1,\ldots,s \\
        & \alpha_i \geq 0,
            & i=1,\ldots,s
\end{align*}
where $C$ is a configurable penalty term.
Then, it assigns label $1$ to $\tilde{p}$ if $\pi^T \tilde{p} \geq \pi_0 $ and $-1$ otherwise. This description is commonly referred to \emph{SVMs with linear kernels}. Non-linear extensions have also been described. In Subsection~\ref{subsec:aff}, we use SVMs to predict affine subspaces.

\subsection{Model Evaluation and Cross-Validation}

To evaluate the performance of a prediction method on a specific training dataset, we employ \emph{$k$-fold cross cross-validation}. Given a training dataset $\mathcal{D} = \{(p^1,h_1),\ldots,(p^s,h_s)\}$, first we partition $\mathcal{D}$ into $k$ disjoint subsets $\mathcal{D}_1,\ldots,\mathcal{D}_k$ of equal size (we assume that $s$ is divisible by $k$) and train $k$ predictors $\phi_1,\ldots,\phi_k$, using the set
$
    \bigcup_{j \in \{1,\ldots,k\} \setminus \{i\}} \mathcal{D}_j
$
as training dataset for predictor $\phi_i$.
The performance of the method on $\mathcal{D}$ is then defined as
\[
    \sum_{i=1}^k \sum_{(p,h) \in \mathcal{D}_i}  \frac{f\left( \phi_i(p), h \right)}{s},
\]
where $f:\{-1,1\}^d \times \{-1,1\}^d \to [0,1]$ is a performance measure. The idea, for each sample $(p,h) \in \mathcal{D}_i$, is to compare the actual label $h$ against the predicted label $\phi_i(p)$, produced by a predictor $\phi_i$ that did not have access to the pair $(p,h)$ during training.
In this work, we use two classical performance measures: \emph{precision} and \emph{recall}, given, respectively, by
\begin{align*}
    f_\text{precision}(\tilde{h}, h) = \frac{\left| \{j : \tilde{h}_j = 1 \land h_j = 1 \} \right|}{\left| \{j : \tilde{h}_j = 1 \} \right|}, \\
    f_\text{recall}(\tilde{h}, h) = \frac{
            \left| \{j : \tilde{h}_j = 1 \land h_j = 1 \} \right|
        }{
            \left| \{j : h_j = 1 \} \right|
        }.
\end{align*}

\section{Setting and Learning Strategies}\label{sec:learning}

In this section, we formalize the setting in which the learning takes place, in addition to our proposed learning strategies. The theoretical setting below was designed to match the realistic conditions that energy market operators face on a daily basis.

Assume that every instance $I(p)$ of SCUC is completely described by a vector of real-valued parameters $p$. These parameters include only the data that is unknown to the market operators the day before the market clearing process takes places. For example, they may include the peak system demand and the production costs for each generator. In this work, we assume that network topology and the total number of generators is known ahead of time, so there is no need to encode them via parameters. We may, therefore, assume that $p \in \mathbb{R}^n$, where $n$ is fixed.
We are given a set of training instances $I(p^1),\ldots,I(p^s)$, which are similar to the instance we expect to solve during actual market clearing.
These training instances can either come directly from historical data, if enough similar instances have been solved in the past, or they can be artificially generated, by sampling the probability distribution of $p$.
Knowledge about this distribution is not an unreasonable assumption, since market operators have been solving SCUC daily for years and have large amounts of accumulated data which can be readily analyzed.

We also assume that we have at our disposal a customized MIP solver which can receive, in addition to the instance $I(p)$ to be solved, a vector of hints. These hints may affect the performance of the MIP solver and will, hopefully, improve its running time. The hints may also affect the quality or optimality of the solutions produced by the MIP solver. Examples of hints may include a list of warm starts, a list variables that should be fixed, or a list of initial constraints to enforce.

The day before market clearing, a \emph{training phase} takes place.
During this phase, our task is to construct a prediction function $\phi$ which maps each possible parameter $p \in \mathbb{R}^n$ into a vector of hints. In order to build this function, we solve the training instances $I(p^1),\ldots,I(p^s)$, capturing any data we may consider useful.
During actual market clearing, when we are under time pressure to solve new instances of the problem, the \emph{test phase} begins.
In this phase, one particular vector of parameters $\tilde{p}$ is given, the instance $I(\tilde{p})$ is constructed, and our goal is to solve $I(\tilde{p})$ very quickly. To do this, the vector of hints $\phi(\tilde{p})$ is computed and the pair $(I(\tilde{p}), \phi(\tilde{p}))$ is given to the MIP solver. The total running time is measured, including the time taken to construct the vector of hints from the parameters, as well as the time taken by the MIP solver to solve the instance.

In the following, we propose three different ML models, which target different challenges of solving SCUC. The first model is focused on predicting violated transmission and N-1 security constraints. The second model focuses on producing initial feasible solutions, which can be used as warm starts. The third model aims at predicting an affine subspace where the solution is likely to lie.

\subsection{Learning Violated Transmission Constraints}
\label{subsec:tr}

As explained in Subection~\ref{subsec:scuc}, one of the most complicating factors when solving SCUC is handling the large number of transmission and security constraints that need to be enforced.
Our first model is designed to predict which transmission constraints should be initially added to the relaxation and which constraints can be safely omitted.

As baseline, we use the iterative contingency screening algorithm presented in \citep{XavierFengWangThimmapuram2018}, which was independently implemented and evaluated at MISO, where it presented better performance than internal methods. The method works by initially solving a relaxed version of SCUC where no transmission or N-1 security constraints are enforced. At each iteration, all violated constraints are identified, but only a small and carefully selected subset is added to the relaxation. The method stops when no further violations are detected, in which case the original problem is solved.

In the following, we modify this baseline method to work with hints from a machine learning predictor.
Let $L$ be the set of transmission lines.
We recall that each transmission constraint takes the form
\begin{align*}
  & -F^c_l \leq \sum_{b \in B} \delta^c_{lb} \left( \sum_{g \in G_b} y_{gt} - d_{bt} \right) \leq F^c_l & \forall l \in L, c \in L \cup \{0\}, t \in T.
\end{align*}
We assume that the customized MIP solver accepts a vector of hints
\begin{align*}
  & h_{l,c,t} \in \{\texttt{ENFORCE}, \texttt{RELAX}\} & \forall l \in L, c \in L \cup \{0\}, t \in T,
\end{align*}
indicating whether the thermal limits of transmission line $l$, under contingency scenario $c$, at time $t$, is likely to be binding in the optimal solution.
Given these hints, the solver proceeds as described in Algorithm~\ref{alg}.
In the first iteration, only the transmission constraints specified by the hints are enforced. For the remaining iterations, the solver follows exactly the same strategy as the baseline method, and keeps adding small subsets of violated constraints until no further violations are found.
If the predictor is successful and correctly predicts a superset of the binding transmission constraints, then only one iteration is necessary for the algorithm to converge. If the predictor fails to identify some necessary constraints, the algorithm may require additional iterations, but still returns the correct solution.
The predictor should still strive to keep a low false positive rate, since the inclusion of too many redundant constraints during the first iteration of the procedure may significantly slow down the MIP optimization time.

\begin{algorithm}
    \footnotesize
    \caption{\footnotesize Security-Constrained Unit Commitment \label{alg}}
    \begin{algorithmic}[1]
        \State Create a relaxation of SCUC without constraints \eqref{c:transmission} or \eqref{c:security}.
        \State Query transmission predictor and receive vector of hints $h$.
        \For {$(l,c,t) \in L \times (L \cup \{0\}) \times T$}
            \State If $h_{l,c,t} = \texttt{ENFORCE}$, add corresponding constraint \eqref{c:transmission} or \eqref{c:security}.
        \EndFor

        \State Solve the current relaxation.
        \For {$(l,c,t) \in L \times (L \cup \{0\}) \times T$}
            \State Let $f^c_{lt} = \sum_{b \in B} \delta^c_{lb} \left( \sum_{g \in G_b} y_{gt} - d_{bt} \right)$.
            \State Let $\gamma^c_{lt} = \max\{- f^c_{lt} - F^c_l, 0, f^c_{lt} - F^c_l\}$.
        \EndFor
        \State Let $\Gamma = \{(l,c,t) \in L \times (L \cup \{0\}) \times T: \gamma^c_{lt} > 0\}$ be set of violated constraints.

        \If {$\Gamma$ is empty}
            \Return current solution.
        \Else
            \State For every $l \in L, t \in T$, keep in $\Gamma$ only the pair $(l, c, t)$ with highest $\gamma^{c}_{lt}$.
            \State For every $t \in T$, keep in $\Gamma$ only the 15 pairs $(l,c,t)$ having the highest $\gamma^c_{lt}$.
            \State For every violation in $\Gamma$, add the corresponding constraint to the relaxation.
            \State \textbf{goto} step 5.
        \EndIf
    \end{algorithmic}
\end{algorithm}

For this prediction task, we propose a simple learning strategy based on $k$-Nearest Neighbors. Let $p^1,\ldots,p^s$ be the training parameters. During training, the instances $I(p^1),\ldots,I(p^s)$ are solved using Algorithm~\ref{alg}, where the initial hints are set to $\texttt{RELAX}$ for all constraints. During the solution of each training instance $I(p^i)$, we record the set of constraints $\mathcal{L}^i \subseteq L \times (L \cup \{0\}) \times T$ that was added to the relaxation by Algorithm~\ref{alg}. During test, given a vector of parameters $\tilde{p}$, the predictor finds a list of the $k$ vectors $p^{i_1},\ldots,p^{i_{k}}$ which are the closest to $\tilde{p}$ and sets $h_{l,c,t}=\texttt{ENFORCE}$ if $(l,c,t)$ appears in at least $p\%$ of the sets $\mathcal{L}^{i_1},\ldots,\mathcal{L}^{i_k}$, where $p$ is a hyperparameter. When $p=50$, this strategy corresponds to the traditional $k$NN algorithm with simple majority voting. Note that the cost of incorrectly identifying a constraint as necessary (false positive) is typically much smaller than the cost of incorrectly identifying a constraint as redundant (false negative): in the former case, the dimension of the problem simply increases a little, while in the latter, an additional iteration may be necessary. For this reason, in our experiments we try to avoid false negatives (at the cost of more false positives) by setting $p=10$. That is, a constraint is added if is necessary for at least 10\% of the nearest training examples.

Note that this $k$NN description generalizes a number of other very simple strategies. For example, when $k=1$, the predictor looks for the most similar instance in the training data. When $k=s$, the strategy is equivalent to counting how many times each transmission constraint was necessary during the solution of the entire training dataset. When $k=1$ and $p=0$, the predictor simply memorizes all constraints that were ever necessary during training. In Subsection~\ref{subsec:eval_tr}, we run experiments for several of these variations. As we discuss later, these simple $k$NN strategies already performed very well in our experiments, which is the reason we do not recommend more elaborated methods.


We note that this $k$NN transmission predictor is also very suitable for online learning, and can be used even when only a very small number of training samples is available. During the solution of the training instances, for example, after a small number of samples has been solved, subsequent samples can already use hints from this predictor to accelerate the solution process.

\subsection{Learning Warm Starts}
\label{subsec:ws}

Another considerably challenging aspect of solving large-scale SCUC instances is obtaining high-quality initial feasible solutions to the problem. Over the last decades, researchers have proposed a number of strong mixed-integer linear programming formulations for SCUC, leading to very strong dual bounds being obtained fairly quickly during the optimization process. For challenging instances, however, we have observed that a significant portion of the running time is usually spent in finding a feasible solution that has an objective value close to this dual bound. These solutions are currently typically found through generic primal heuristics included in the MIP solver, such as feasibility pump \citep{FischettiGloverLodi2005} or large neighborhood search \citep{Shaw1998}.
Modern MIP solvers also allow users to specify user-constructed solutions, which can be used by the solver to accelerate various parts of the branch-and-cut process. Very importantly, the solutions provided by the user do not need to be complete, or even feasible. Modern solvers can solve restricted subproblems to find values of missing variables, or even repair infeasible solutions using various heuristics. The strategy we present in this subsection takes advantage of these capabilities, instead of trying to replace them. Our aim is to produce, based on a large number of previous optimal solution, a (partial) solution which works well as warm start.

For this prediction task we propose, once again, a strategy based on $k$-Nearest Neighbors.
Let $p^1,\ldots,p^s$ be the training parameters and let $(x^1,y^1),\ldots, (x^s,y^s)$ be their optimal solutions. In principle, the entire list of solutions $(x^1,y^1),\ldots,(x^s,y^s)$ could be provided as warm starts to the solver; however, processing each solution requires a non-negligible amount of time, as we show in Subsection~\ref{subsec:eval_ws}. We propose instead the usage of $k$NN to construct a single (partial and possibly infeasible) solution, which should be repaired by the solver. The idea is to set the values only for the variables which we can predict with high confidence, and to let the MIP solver determine the values for all the remaining variables. To maximize the likelihood of our warm starts being useful, we focus on predicting values for the binary variables $x$; all continuous variables $y$ are determined by the solver. Given a test parameter $\tilde{p}$, first we find the $k$ solutions $(x^{i_1},y^{i_1}),\ldots, (x^{i_k},y^{i_k})$ corresponding to the $k$ nearest training instances $I(p^{i_1}),\ldots,I(p^{i_k})$. Then, for each variable $x_{gt}$, we verify the average $\tilde{x}_{gt}$ of the values $x^{i_1}_{gt}, \ldots, x^{i_k}_{gt}$. If  $\tilde{x}_{gt} > p$, where $p \in [0.5, 1.0]$ is a hyperparameter, then the predictor sets $x_{g,t}$ to one in the warm start solution. If $\tilde{x}_{gt} \leq (1-p)$, then the value is set to zero. In any other case, the value is left undefined.

When $p=0.5$, this strategy corresponds to the traditional $k$NN method, where all variables are defined, even if the neighboring solutions significantly disagree.
As $p$ increases, a higher threshold for a consensus is required, and variables for which consensus is not reached are left undefined.
More complete warm starts can be processed significantly faster by the solver, but have higher likelihood of being sub-optimal or irreparably infeasible.
In Subsection~\ref{subsec:eval_ws}, we evaluate several choices of $p$ to quantify this trade-off. In our computational experiments, $p=0.9$ provided the best results.


Like the transmission predictor presented previously, this warm-start predictor is also very suitable for online learning, and can be used with a very small number of training samples.

\subsection{Learning Affine Subspaces}
\label{subsec:aff}

Through their experience, market operators have naturally learned a number of patterns in the optimal solutions of their own power systems. For example, they know that certain power generating units (known as \emph{base units}) will almost certainly be committed throughout the day, while other units (known as \emph{peak units}) typically come online only during peak demand. These and many other intuitive patterns, however, are completely lost when SCUC is translated into a mathematical problem. Given the fixed characteristics of a particular power system and a realistic parameter distribution, new constraints could probably be safely added to the MIP formulation of SCUC to restrict its solution space without affecting its set of optimal solutions.
In this subsection, we develop a ML model for finding such constraints automatically, based on statistical data.

More precisely, let $\tilde{p}$ be a vector of parameters. Our goal is to develop a model which predicts a list of hyperplanes $(h^1,h^1_0),\ldots,(h^k,h^k_0)$ such that, with very high likelihood, the optimal solution $(x,y)$ of $I(\tilde{p})$ satisfies $\langle h^i, x \rangle = h^i_0$, for $i=1,\ldots,k$.
In this work, we restrict ourselves to hyperplanes that can be written as
\begin{align}
  x_{gt} & = 0, \label{eq:fixzero} \\
  x_{gt} & = 1, \label{eq:fixone} \\
  x_{gt} & = x_{g,t+1}, \label{eq:fixnext}
\end{align}
where $g \in G$ and $t \in T$. In other words, the predictor tries to determine whether each variable $x_{gt}$ should be fixed to zero, to one, or to the next variable $x_{g,t+1}$. Furthermore, to prevent conflicting assignments the predictor makes at most one of these three recommendations for each variable $x_{g,t}$.

Let $\mathcal{H}$ be the set of all hyperplanes considered for inclusion.
During the training phase, our goal is to construct, for each hyperplane $(h,h_0) \in \mathcal{H}$,  a prediction function $\phi_{(h,h_0)} : \mathbb{R}^n \to \{\texttt{ADD}, \texttt{SKIP}\}$, which receives a vector of parameters $\tilde{p} \in \mathbb{R}^n$ and returns the label \texttt{ADD} to indicate that the equality constraint $\langle h, x \rangle = h_0$ is very likely to be satisfied at the optimal solution, and \texttt{SKIP} otherwise.
Given these hints and a vector of parameteres $p$, our customized MIP solver simply adds to the relaxation all the equality constraints $\langle h, x \rangle = h_0$ such that $\phi_{(h,h_0)}(p) = \texttt{ADD}$. These constraints are added at the very beginning of the optimization process and kept until the final solution is obtained; that is, they are never removed. While, in principle, adding these constraints could lead to sub-optimal solutions, in our computational experiments we have observed that, by using a reasonably large number of training samples and by carefully constructing high-precision predictors, as described in detail below, this strategy does not lead to any noticeable degradation in solution quality, while providing significant speedups.

Next we describe how can such high-precision predictors can be constructed in practice. Let $(h,h_0) \in \mathcal{H}$. Furthermore, let $I(p^1),\ldots,I(p^s)$ be the training instances and let $(x^1,y^1),\ldots,(x^s,y^s)$ be their respective optimal solutions. For every $i \in \{1,\ldots,s\}$, let $z^i \in \{0,1\}$ be the label indicating whether solution $(x^i,y^i)$ lies in the hyperplane $(h,h_0)$ or not. That is,
$z^i = 1$ if $\langle h, x^i \rangle = h_0$, and $z^i = 0$ otherwise.
We also denote by $\bar{z}$ the average value of the variables $z^1, \ldots, z^s$.

The proposed algorithm for constructing the function $\phi_{(h,h_0)}$ is the following.
First, if $\bar{z}$ is very close to one, then the predictor always suggests adding this hyperplane to the formulation. More precisely, if $\bar{z} \geq z^{\text{FIX}}$, where $z^{\text{FIX}}$ is a fixed hyper-parameter, then the prediction function $\phi_{(h,h_0)}$ returns \texttt{ADD} for every input.
Next, the the algorithm verifies whether the labels are sufficiently balanced to reliably perform supervised training. More precisely, if $\bar{z} \notin [ z^{\text{MIN}}, z^{\text{MAX}} ]$, where $z^{\text{MIN}}$ and $z^{\text{MAX}}$ are hyper-parameters, then the prediction function $\phi_{(h,h_0)}$ always returns \texttt{SKIP}.
Finally, if $\bar{z} \in [ z^{\text{MIN}}, z^{\text{MAX}} ]$, then the algorithm constructs a binary classifier and evaluates its precision and recall on the training set using $k$-fold cross validation.
If the binary classifier proves to have sufficiently high precision and recall, the prediction function $\phi_{(h,h_0)}$ returns its predictions during the test phase. Otherwise, the algorithm discards the binary classifier, and the function $\phi_{(h,h_0)}$ returns \texttt{SKIP} for every input. The minimum acceptable recall is given by a fixed hyperparameter $\alpha^{\text{R}}$. The minimum acceptable precision is computed by the expression
\[
    \max\{\bar{z}, 1 - \bar{z}\} \cdot (1 - \alpha^{\text{P}}) + \alpha^{\text{P}},
\]
where $\alpha^{\text{P}}$ is a hyper-parameter. Intuitively, the algorithm only accepts the trained binary classifier if it significantly outperforms a dummy classifier which always returns the same label for every input. When $\alpha^{\text{P}}=0$, the trained binary classifier is acceptable as long as it has same precision as the dummy classifier. When $\alpha^{\text{P}}=1$, the classifier is only acceptable if it has perfect precision.

Table~\ref{table:hyp} shows the precise hyper-parameters used in our experiments, for each type of hyperplane. Due to the high dimensionality of the parameters, we do not train the binary classifiers directly the original vector of parameters $p$. Instead, we propose using as training features only the following subset of parameters: the peak system load, the hourly system loads, the average production cost of generator $g$ and the average production costs of the remaining generators. We propose the usage of support vector machines (with linear kernels), since, in our experiments, these models performed better than logistic regression and random decision forests models, while remaining computationally friendly. Neural networks were not considered given the small number of training samples available.

\begin{table}
  \TABLE
  {Hyper-parameters for affine subspace predictor. \label{table:hyp}}
  {
    \centering \small
    \setlength{\tabcolsep}{0.5em}
    \begin{tabular}{lrrrrr}
    \toprule
    Hyperplane
    & $z^{\text{FIX}}$
    & $z^{\text{MIN}}$
    & $z^{\text{MAX}}$
    & $\alpha^{\text{R}}$
    & $\alpha^{\text{P}}$
    \\
    \midrule
    $x_{gt} = 0$         & 1.000 & 0.250 & 0.750 & 0.90 & 0.90 \\
    $x_{gt} = 1$         & 1.000 & 0.250 & 0.750 & 0.75 & 0.75 \\
    $x_{gt} = x_{g,t+1}$ & 0.975 & 0.025 & 0.975 & 0.50 & 0.75 \\
    \midrule
  \end{tabular}
  }
  {}
\end{table}

Unlike the previous two predictors, the affine subspace predictor described in this subsection requires a substantial number of samples in order to give reliable predictions, making it less suitable for online learning.

\section{Computational Experiments}\label{sec:comp}

In this section we evaluate the practical impact of the ML models introduced in Section~\ref{sec:learning} by performing extensive computational experiments on a diverse set of large-scale and realistic SCUC instances.
In Subsections~\ref{subsec:instances} and \ref{subsec:samples}, we describe the our implementation of the proposed methods, the computational environment, and the instances used for benchmark. In Subsections~\ref{subsec:eval_tr}, \ref{subsec:eval_ws} and \ref{subsec:eval_aff} we evaluate, respectively, the performance of the three predictors proposed in Section~\ref{sec:learning}, in addition to several variations. In Subsection~\ref{subsec:ood}, we evaluate the out-of-distribution performance of these predictors, to measure their robustness against moderate dataset shift.

\subsection{Computational Environment and Instances}
\label{subsec:instances}

The three proposed machine-learning models described in Section~\ref{sec:learning} were implemented in Julia 1.2 \citep{julialang} and \texttt{scikit-learn} \citep{scikitlearn}.
Package ScikitLearn.jl \citep{scikitlearnjl} was used to access \texttt{scikit-learn} predictors from Julia, while package DataFrames.jl \citep{dataframesjl} was used for tabular data manipulation.
IBM ILOG CPLEX 12.8.0 \citep{cplex} was used as the MIP solver. The code responsible for loading the instances, constructing the MIP, querying the machine learning models and translating the given hints into instructions for CPLEX was also written in Julia, using JuMP \citep{JuMP} as modeling language. Solver features not available through JuMP were accessed through CPLEX's C API. Training instances were generated and solved at the Bebop cluster at Argonne National Laboratory (Intel Xeon E5-2695v4 2.10GHz, 18 cores, 36 threads, 128GB DDR4). During the training phase, multiple training instances were solved in parallel at a time, on multiple nodes. To accurately measure the running times during the test phase, a dedicated node at the same cluster was used.
CPLEX was configured to use a relative MIP gap tolerance of 0.1\% during test and 0.01\% during training.

\begin{table}
  \TABLE
  {Size of selected instances. \label{table:size}}
  {
    \centering \small
    \begin{tabular}{lrrr}
    \toprule
    Instance
    & \hspace{1em} Buses
    & \hspace{1em} Units
    & \hspace{1em} Lines
    \\
    \midrule
    \texttt{case1888rte} & 1,888 &   297 & 2,531 \\
    \texttt{case1951rte} & 1,951 &   391 & 2,596 \\
    \texttt{case2848rte} & 2,848 &   547 & 3,776 \\
    \texttt{case3012wp}  & 3,012 &   502 & 3,572 \\
    \texttt{case3375wp}  & 3,374 &   596 & 4,161 \\
    \texttt{case6468rte} & 6,468 & 1,295 & 9,000 \\
    \texttt{case6470rte} & 6,470 & 1,330 & 9,005 \\
    \texttt{case6495rte} & 6,495 & 1,372 & 9,019 \\
    \texttt{case6515rte} & 6,515 & 1,388 & 9,037 \\
    \bottomrule
    \end{tabular}
  }
  {}
\end{table}

A total of nine realistic instances from MATPOWER \citep{matpower} were selected to evaluate the effectiveness of the method. These instances correspond to realistic, large-scale European test systems. Table~\ref{table:size} presents their main characteristics, including the number of buses, units and transmission lines. Some generator data necessary for SCUC was missing in these instances, and was artificially generated based on publicly available data from \cite{PjmWebsite} and \cite{MisoWebsite}.

\subsection{Training and Test Samples}
\label{subsec:samples}

During training, 300 variations for each of the nine original instances were generated and solved. We considered four sources of uncertainty: i) uncertain production and startup costs; ii) uncertain geographical load distribution; iii) uncertain peak system load; and iv) uncertain temporal load profile.
The precise randomization scheme is described below. The specific parameters were chosen based on our analysis of publicly available bid and hourly demand data from \cite{PjmWebsite}, corresponding to the month of January, 2017.

\begin{figure}
    \FIGURE
    {\includegraphics[width=0.75\textwidth]{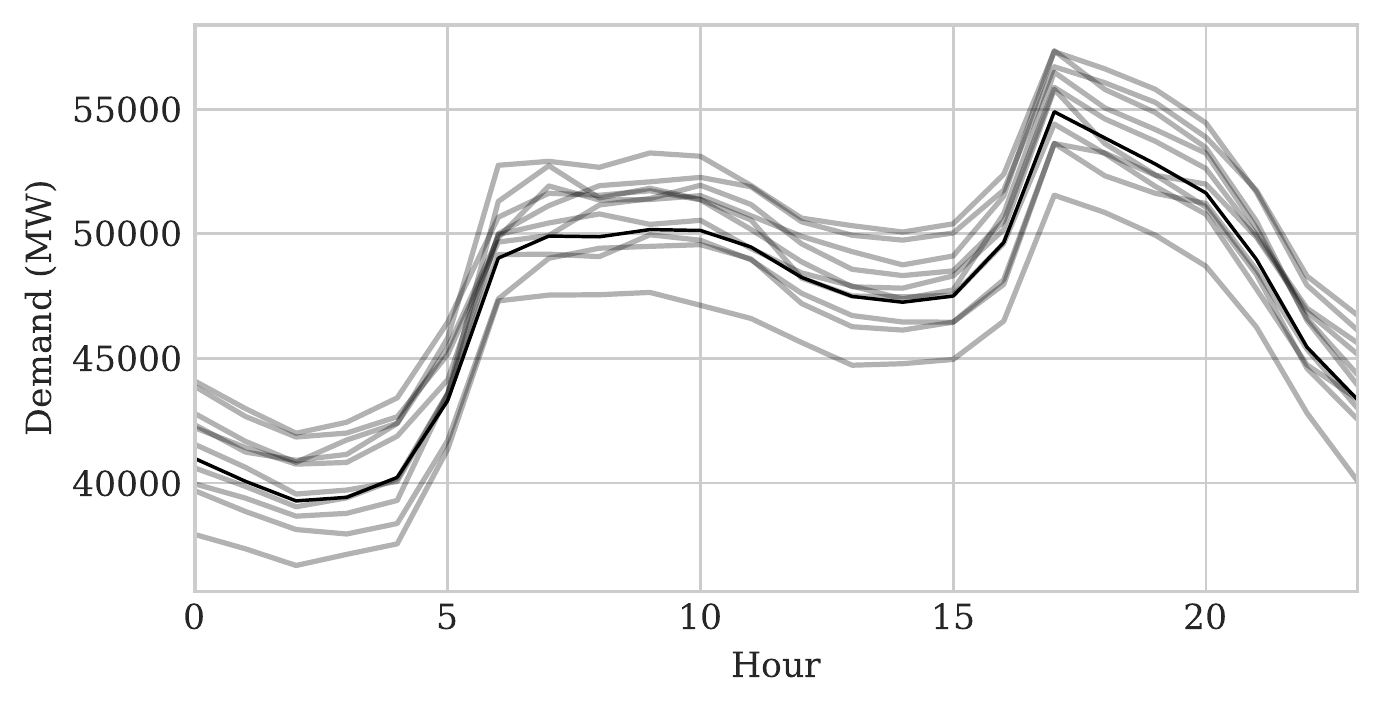}}
    {Sample of artificially generated load profiles. \label{fig:load-profile}}
    {}
\end{figure}

\begin{enumerate}[(i)]
    \item \textbf{Production and startup costs}. In the original instances, the production cost for each generator $g \in G$ is modelled as a convex piecewise-linear function described by the parameters $c^0_g$, the cost of producing the minimum amount of power, and $c^1_g, \ldots, c^k_g$, the costs of producing each additional MW of power within the piecewise interval $k$. In addition, each generator has a startup cost $c^s_g$.
    In our data analysis, we observed that the daily changes in bid prices rarely exceeded $\pm 5\%$. Therefore, random numbers $\alpha_g$ were independently drawn from the uniform distribution in the interval $[0.95, 1.05]$, for each generator $g \in G$, and the costs of $g$ were linearly scaled by this factor. That is, the costs for generator $g$ were set to $\alpha_g c^0_g, \alpha_g c^1_g, \ldots, \alpha_g c^k_g, \alpha_g c^s_g$.

    \item \textbf{Geographical load distribution}. In the original instances, each bus $b \in B$ is responsible for a certain percentage $d_b$ of the total system load. To generate variations of these parameters, random numbers $\beta_b$ were independently drawn from the uniform distribution in the interval $[0.90, 1.10]$.
    The percentages $d_b$ were then linearly scaled by the $\beta_b$ factors and later normalized. More precisely, the demand for each bus $b \in B$ was set to $\frac{\beta_b d_b}{\sum_{i \in b} \beta_i}$.

    \item \textbf{Peak system load and temporal load profile}.
    For simplicity, assume $T=\{1,\ldots,24\}$. Let $D_1,\ldots,D_{24}$ denote the system-wide load during each hour of operation, and let $\gamma_t = \frac{D_{t+1}}{D_t}$ for $t=1,\ldots,23$, be the hourly variation in system load.
    In order to generate realistic distribution for these parameters, we analyzed hourly demand data from \cite{PjmWebsite}. More specifically, we computed the mean $\mu_t$ and variance $\sigma^2_t$ of each parameter $\gamma_t$, for $t=1,\ldots,23$. To generate instance variations, random numbers $\gamma'_t$ were then independently drawn from the Gaussian distribution $\mathcal{N}(\mu_t, \sigma^2_t)$. Note that the $\gamma'$ parameters only specify the variation in system load from hour to hour, and are not sufficient to determine $D_1,\ldots,D_{24}$.
    Therefore, in addition to these parameters, a random number $\rho$, corresponding to the peak system load $\max\{D_1,\ldots,D_{24}\}$, was also generated. In the original instances, the peak system load is always 60\% of the total capacity. Based on our data analysis, we observed that the actual peak load rarely deviates more than $\pm 7.5\%$ from the day-ahead forecast. Therefore, to generate instance variations, $\rho$ was sampled from the uniform distribution in the interval $[0.6 \times 0.925C, 0.6 \times 1.075 C]$, where $C$ is the total capacity. Note that the $\rho$ and $\gamma'$ parameters are now sufficient to construct $D_1,\ldots,D_{24}$. Figure~\ref{fig:load-profile} shows a sample of some artificially generated load profiles.
\end{enumerate}
In this work, we do not consider changes to the network topology, since this information is assumed to be know by the operators the day before the market clearing takes places. If the network topology is uncertain, an additional parameter $p_l \in \{0,1\}$ can be added for each transmission line $l \in L$, indicating whether the line is available or not.

During the test phase, 50 additional variations of each original instance were generated. In Subsection~\ref{subsec:eval_tr}, \ref{subsec:eval_ws} and \ref{subsec:eval_aff}, the test samples were generated using the same randomization scheme outline above, but with a different random seed.
In Subsection~\ref{subsec:ood}, as described in more detailed in that section, the test samples were generated based on a similar, but different distribution.

\subsection{Evaluation of Transmission Predictor}
\label{subsec:eval_tr}

We start by evaluating the performance of different strategies for predicting necessary transmission and N-1 security constraints. For this prediction task, we compared seven different strategies. Strategy \texttt{zero} corresponds to the case where no machine learning is used. In strategy \texttt{tr:nearest}, we add to the initial relaxation all constraints that were necessary during the solution of the most similar training variation. In \texttt{tr:all}, we add all constraints that were ever necessary during training.
Strategies \texttt{tr:knn:$k$} correspond to the method outlined in Subsection~\ref{subsec:tr}, where $k \in \{10,50,300\}$ is the hyperparameter describing the number of neighbors. We recall that a constraint is initially added to the relaxation it was necessary for at least 10\% of the nearest neighbors.
Since there are only 300 variations for each instance, strategy \texttt{tr:knn:300} is equivalent to simply counting, over all variations, how many times each constraint was necessary.
Finally, strategy \texttt{tr:perf} shows the performance of the theoretically perfect predictor, which can predict the exact set of constraints found by the iterative method, with perfect accuracy, under zero running time. To implement this predictor, we gave it access to the optimal solution produced by the \texttt{zero} strategy.
For each strategy, Table~\ref{table:trAccuracy} shows the average wallclock running time (in seconds), the average number of constraints added to the relaxation per time period, as well as the average number of iterations required by Algorithm~\ref{alg} to converge.
Table~\ref{table:trSpeed} shows a more detailed breakdown of the average running times per instance.
\begin{table}
  \TABLE
  {Impact of different transmission predictors on number of constraints added per time period, number of iterations and running time. \label{table:trAccuracy}}
  {
    \setlength{\tabcolsep}{0.5em}
    \begin{tabular}{lrrrrrrr}
        \toprule
        Method &
        Iterations &
        Violations &
        Time (s) &
        Speedup
        \\
        \midrule
        \begin{filecontents*}{tables/tr_accuracy.csv}
algorithm,wallclockTime,iterations,violations,speedupZero
zero,226.67,4.92,11.10,1.00
tr:nearest,110.62,1.63,12.14,2.05
tr:all,91.27,1.01,25.12,2.48
tr:knn:10,84.53,1.06,16.36,2.68
tr:knn:50,82.85,1.02,15.74,2.74
tr:knn:300,80.36,1.02,15.92,2.82
tr:perf,80.94,1.00,11.10,2.80
\end{filecontents*}
\csvreader[
                head to column names,
                late after line={\\},
                late after last line={\\},
                filter expr={
                    test{\ifnumless{\thecsvinputline}{11}}
                }
            ]{tables/tr_accuracy.csv}{}{%
                \texttt{\algorithm} &
                \iterations &
                \violations &
                \wallclockTime &
                {\speedupZero}x &
            }%
        \bottomrule
    \end{tabular}
  }
  {}
\end{table}
\begin{table}
  \TABLE
  {Impact of different transmission predictors on total wallclock time. \label{table:trSpeed}}
  {
    \setlength{\tabcolsep}{0.5em}
    \begin{tabular}{lrrrrrrr}
        \toprule
        & \multicolumn{7}{c}{Wallclock Time (s)}
        \\
        \cmidrule(r){2-8}
        Instance
        & \texttt{zero}
        & \texttt{tr:nearest}
        & \texttt{tr:all}
        & \texttt{tr:knn:10}
        & \texttt{tr:knn:50}
        & \texttt{tr:knn:300}
        & \texttt{tr:perf}
        \\
        \midrule
        \begin{filecontents*}{tables/tr_speed.csv}
instance,wallclockTime:tr:all,wallclockTime:tr:knn:10,wallclockTime:tr:knn:300,wallclockTime:tr:knn:50,wallclockTime:tr:nearest,wallclockTime:tr:perf,wallclockTime:zero
case1888rte,10.11,9.22,9.26,9.32,14.80,9.52,35.68
case1951rte,15.34,13.00,12.68,13.05,15.69,12.51,26.07
case2848rte,18.07,17.59,17.16,17.28,25.53,18.14,39.38
case3012wp,30.75,29.22,27.96,28.71,43.95,30.21,66.11
case3375wp,82.85,69.65,64.17,62.98,84.26,65.68,124.05
case6468rte,169.87,164.18,154.04,158.82,233.31,142.21,477.48
case6470rte,98.55,97.41,89.54,89.37,122.12,97.84,290.55
case6495rte,190.93,157.92,154.92,157.37,220.54,158.64,540.38
case6515rte,205.00,202.60,193.56,208.76,235.40,193.71,440.37
Average,91.27,84.53,80.36,82.85,110.62,80.94,226.67
\end{filecontents*}
\csvreader[
                head to column names,
                late after line={\\},
                late after last line={\\},
                filter expr={
                    test{\ifnumless{\thecsvinputline}{11}}
                }
            ]{tables/tr_speed.csv}{
                1=\instance,
                2=\WTrAll,
                3=\wTrKnnTen,
                4=\wTrKnnThreeHundred,
                5=\wTrKnnFifty,
                6=\wTrNearest,
                7=\wTrPerf,
                8=\wTrZero,
            }{%
                \texttt{\instance} &
                \wTrZero &
                \wTrNearest &
                \WTrAll &
                \wTrKnnTen &
                \wTrKnnFifty &
                \wTrKnnThreeHundred &
                \wTrPerf
            }%
        \midrule
        \csvreader[
                head to column names,
                late after line={\\},
                late after last line={\\},
                filter expr={
                    test{\ifnumgreater{\thecsvinputline}{10}}
                }
            ]{tables/tr_speed.csv}{
                1=\instance,
                2=\WTrAll,
                3=\wTrKnnTen,
                4=\wTrKnnThreeHundred,
                5=\wTrKnnFifty,
                6=\wTrNearest,
                7=\wTrPerf,
                8=\wTrZero,
            }{%
                \texttt{\instance} &
                \wTrZero &
                \wTrNearest &
                \WTrAll &
                \wTrKnnTen &
                \wTrKnnFifty &
                \wTrKnnThreeHundred &
                \wTrPerf
            }
        \bottomrule
    \end{tabular}
  }
  {}
\end{table}

Tables~\ref{table:trAccuracy} and \ref{table:trSpeed} show us that the prediction of necessary transmission and N-1 security constraints is a relatively easy prediction task, which, nevertheless, has very high impact on optimization time.
Previous research has already shown that the set of necessary transmission constraints remains relatively unchanged even under significant load variations \citep{RoaldMolzahn2019}. Our experiments here support these findings, and suggest this is true also for N-1 security constraints, and even under significant variations in temporal load profiles and production costs.
All strategies, even the most simple ones, were quite effective, and significantly improved running times, with speedups ranging from 2.0x to 2.8x.
Strategy {\tt tr:nearest} added the smallest number of constraints to the relaxation, however it often failed to identify some necessary constraints, requiring additional iterations. Strategy {\tt tr:all} was the most effective at reducing the number of iterations, however it was somewhat over-conservative, and added significantly more constraints to the relaxation than the alternatives. The {\tt tr:knn} strategies achieved a better balance between low number of constraints and low number of iterations. Among all practical strategies studied, {\tt tr:knn:300} presented the best average running times by a small margin. Overall, all methods performed very close to {\tt tr:perf}. We note that the average running time of {\tt tr:perf} was slightly higher than {\tt tr:knn:300}, due to the natural performance variability of solving MIPs \citep{LodiTramontani2013}, but we do not consider this difference significant.

Table~\ref{table:trAccuracy} also shows the large negative performance impact of transmission constraints in SCUC. Comparing {\tt tr:all} and {\tt tr:perf}, we see that the inclusion of just 14 redundant constraints per time period (a very small number when compared to the total number of constraints in the problem, which ranges from 55,000 to 230,000 after preprocessing) is already sufficient to induce noticeable slowdown in MIP performance. This is also verified by comparing the running times of instances {\tt case6468rte} and {\tt case6470rte}. Although these instances have roughly the same number of buses, transmission lines and generators, the transmission system in {\tt case6468rte} is significantly more congested, making it almost twice as hard to solve.
The negative impact of transmission constraints in congested systems was also observed by \cite{ChenCastoWangWangWangWan2016}.

\subsection{Evaluation of Warm Start Predictor}
\label{subsec:eval_ws}

Next we focus on the performance of the warm start predictor described in Subsection~\ref{subsec:ws}. For each strategy, Table~\ref{table:wsAccuracy} shows how often was the strategy successful at producing at least one valid warm start, as well as the average wallclock running time (in seconds) to solve the problem. The table also shows how good, in relative terms, the best warm start was, when compared to the optimal solution eventually found by the solver. For example, a gap of 0.25\% in the table indicates that the warm start had objective value 0.25\% higher than the optimal objective value eventually obtained by the solver. A gap of 0\% indicates that the warm start turned out to be optimal. All warm starts were provided to CPLEX with effort level {\tt CPX\_MIPSTART\_REPAIR}, which instructs the solver to attempt to repair the warm start in case it is infeasible, and to solve a sub-MIP if the the values of some integer variables are missing. All other CPLEX parameters related to warm starts were left as default. Table~\ref{table:wsSpeed} shows a more detailed breakdown of running times.

For this prediction task, we compared six different strategies. Strategy {\tt tr:knn:300}, as introduced in Subsection~\ref{subsec:eval_tr}, is included as a baseline, and only predicts transmission constraints, not warm starts. All other strategies presented in this subsection are built on top of {\tt tr:knn:300}; that is, prediction of transmission constraints is performed first, using {\tt tr:knn:300}, then prediction of warm starts is performed. In strategies {\tt ws:collect:$n$}, we provide to the solver $n \in \{5,15\}$ warm starts, containing the optimal solutions of the $n$ nearest training variations. To maximize the likelihood of the warm starts being useful, we only provide values for the integer variables. Strategies {\tt ws:knn:$k$:$p$}, for $k=50$ and $p \in \{50,75,90\}$ correspond to the strategy described in Subsection~\ref{subsec:ws}. To recall, this strategy collects the solutions to the $k$ nearest training variations and constructs a single warm start, where the value of each binary variable is provided to the solver only if the collected solutions reach a $p\%$ consensus on that variable. Strategy {\tt ws:knn:50} corresponds to the traditional $k$-Nearest Neighbors algorithm, using simple majority vote. Finally, strategy {\tt ws:perf} shows the performance of the theoretically perfect predictor, which constructs a single warm start containing the optimal values for all binary variables, under zero running time. To implement this predictor we gave it access to the optimal solution produced by {\tt tr:knn:300}.
\begin{table}
  \TABLE
  {Performance of different warm-start predictors. \label{table:wsAccuracy}}
  {
    \setlength{\tabcolsep}{0.5em}
    \begin{tabular}{lrrrrrrr}
        \toprule
        Method &
        Success (\%) &
        Gap (\%) &
        Time (s) &
        Speedup
        \\
        \midrule
        \begin{filecontents*}{tables/ws_accuracy.csv}
algorithm,wallclockTime,mipstartSuccess,mipstartGap,speedupZero,speedupTr
tr:knn:300,80.36,0.00,---,2.82,1.00
ws:collect:5,64.53,92.22,0.22,3.51,1.25
ws:collect:15,76.37,97.41,0.18,2.97,1.05
ws:knn:50:50,80.10,0.00,---,2.83,1.00
ws:knn:50:75,72.98,46.67,0.08,3.11,1.10
ws:knn:50:90,52.80,100.00,0.00,4.29,1.52
ws:perf,46.06,100.00,0.00,4.92,1.74
\end{filecontents*}
\csvreader[
                head to column names,
                late after line={\\},
                late after last line={\\},
                filter expr={
                    test{\ifnumless{\thecsvinputline}{11}}
                }
            ]{tables/ws_accuracy.csv}{}{%
                \texttt{\algorithm} &
                \mipstartSuccess &
                \mipstartGap &
                \wallclockTime &
                {\speedupZero}x &
            }%
        \bottomrule
    \end{tabular}
  }
  {}
\end{table}
\begin{table}
  \TABLE
  {Impact of different warm-start predictors on total wallclock time. \label{table:wsSpeed}}
  {
    \setlength{\tabcolsep}{0.5em}
    \begin{tabular}{lrrrrrrr}
        \toprule
        & \multicolumn{7}{c}{Wallclock Time (s)}
        \\
        \cmidrule(r){2-8}
        Instance
        & \small\tt tr:knn:50
        & \small\tt ws:collect:5
        & \small\tt ws:collect:15
        & \small\tt ws:knn:50:50
        & \small\tt ws:knn:50:75
        & \small\tt ws:knn:50:90
        & \small\tt ws:perf
        \\
        \midrule
        \begin{filecontents*}{tables/ws_speed.csv}
instance,wallclockTime:tr:knn:300,wallclockTime:ws:collect:15,wallclockTime:ws:collect:5,wallclockTime:ws:knn:50:50,wallclockTime:ws:knn:50:75,wallclockTime:ws:knn:50:90,wallclockTime:ws:perf,speedupTr:tr:knn:300,speedupTr:ws:collect:15,speedupTr:ws:collect:5,speedupTr:ws:knn:50:50,speedupTr:ws:knn:50:75,speedupTr:ws:knn:50:90,speedupTr:ws:perf
case1888rte,9.26,12.37,10.13,9.38,9.20,9.15,9.04,1.00,0.75,0.91,0.99,1.01,1.01,1.02
case1951rte,12.68,15.63,13.07,12.24,11.96,12.50,11.47,1.00,0.81,0.97,1.04,1.06,1.01,1.10
case2848rte,17.16,23.75,20.15,17.72,18.46,18.75,17.04,1.00,0.72,0.85,0.97,0.93,0.92,1.01
case3012wp,27.96,28.85,24.35,29.97,27.30,21.48,18.70,1.00,0.97,1.15,0.93,1.02,1.30,1.49
case3375wp,64.17,60.75,60.40,64.80,54.82,42.24,35.60,1.00,1.06,1.06,0.99,1.17,1.52,1.80
case6468rte,154.04,130.54,102.13,140.61,139.80,99.42,80.38,1.00,1.18,1.51,1.10,1.10,1.55,1.92
case6470rte,89.54,97.29,78.90,97.27,86.49,63.10,58.65,1.00,0.92,1.13,0.92,1.04,1.42,1.53
case6495rte,154.92,145.94,120.95,156.95,153.16,107.12,85.58,1.00,1.06,1.28,0.99,1.01,1.45,1.81
case6515rte,193.56,172.25,150.73,191.94,155.63,101.47,98.06,1.00,1.12,1.28,1.01,1.24,1.91,1.97
Average,80.36,76.37,64.53,80.10,72.98,52.80,46.06,1.00,1.05,1.25,1.00,1.10,1.52,1.74
\end{filecontents*}
\csvreader[
                head to column names,
                late after line={\\},
                late after last line={\\},
                filter expr={
                    test{\ifnumless{\thecsvinputline}{11}}
                }
            ]{tables/ws_speed.csv}{
                1=\instance,
                2=\wTrKnnFifty, 
                3=\wWsCollectFifteen, 
                4=\wWsCollectFive, 
                5=\wWsKnnFiftyFifty, 
                6=\wWsKnnFiftySeventyFive, 
                7=\wWsKnnFiftyNinety, 
                8=\wWsPerf, 
            }{%
                \instance &
                \wTrKnnFifty &
                \wWsCollectFive &
                \wWsCollectFifteen &
                \wWsKnnFiftyFifty &
                \wWsKnnFiftySeventyFive &
                \wWsKnnFiftyNinety &
                \wWsPerf
            }%
        \midrule
        \csvreader[
                head to column names,
                late after line={\\},
                late after last line={\\},
                filter expr={
                    test{\ifnumgreater{\thecsvinputline}{10}}
                }
            ]{tables/ws_speed.csv}{
                1=\instance,
                2=\wTrKnnFifty, 
                3=\wWsCollectFifteen, 
                4=\wWsCollectFive, 
                5=\wWsKnnFiftyFifty, 
                6=\wWsKnnFiftySeventyFive, 
                7=\wWsKnnFiftyNinety, 
                8=\wWsPerf, 
            }{%
                \bf \instance &
                \bf \wTrKnnFifty &
                \bf \wWsCollectFive &
                \bf \wWsCollectFifteen &
                \bf \wWsKnnFiftyFifty &
                \bf \wWsKnnFiftySeventyFive &
                \bf \wWsKnnFiftyNinety &
                \bf \wWsPerf
            }
        \bottomrule
    \end{tabular}
  }
  {}
\end{table}

In contrast to predicting necessary transmission constraints, predicting warm starts proved to be a significantly harder machine learning task. As shown in Table~\ref{table:wsSpeed}, while strategy {\tt ws:collect:5} was somewhat useful for instances with more than 6,000 buses, it was not useful for all the remaining instances. While its achieved speedup of 3.5x was higher than the baseline 2.8x, it was still well below the theoretical maximum of 4.9x. Table~\ref{table:wsAccuracy} shows that, while this strategy produced valid warm starts for 92\% of the problem variations, the warm starts were of relatively low quality, with an average gap of 0.22\%. The table also shows that increasing the number of added warm starts provided to the solver did little to fix this issue. The average gap of strategy {\tt ws:collect:15} was only slightly better, at 0.18\%. The time necessary to process these additional 10 warm starts, however, was significant, and almost completely invalidated any potential benefits provided by them.
Combining previous solutions into a complete warm start also proved challenging. Strategy {\tt ws:knn:50:50}, or simple $k$NN with majority voting, proved ineffective for instances of all sizes. The main issue, in this case, was that none of the combined warm starts produced were feasible, or could even be repaired by CPLEX.
This difficulty in generating complete valid solutions for large-scale SCUC instances was also experienced in previous research efforts trying to solve this problem through machine learning methods alone.

We only obtained more significant success in warm start prediction when we shifted some of the computational burden back to the MIP solver. Although {\tt ws:knn:50:90} does not predict values for all binary variables, we observed that the MIP solver had no difficulty in repairing these missing values. Overall, {\tt ws:knn:50:90} obtained an average speedup of 4.3x, which is much closer to the theoretical maximum of 4.9x than previous strategies. This strategy was very consistent, and produced valid (partial) warm starts for all instance variations. More importantly, after being repaired by CPLEX, these warm starts turned out to be optimal in all cases.

We also experimented with other $k$NN variations, trying to produce more complete warm starts. Overall, the results were negative. We present strategy {\tt ws:knn:50:75} as an illustration. Although the (partial) warm starts provided by this strategy were of high-enough quality (when they could be repaired by CPLEX), the repairing procedure failed too often, hurting the overall performance. Overall, strategy {\tt ws:knn:50:75} was not any better than {\tt ws:knn:50:90}.

\subsection{Evaluation of Affine Subspace Predictor}
\label{subsec:eval_aff}

Finally, we evaluate the performance of different affine subspace predictors. For this prediction task, we compared 6 different methods. Method {\tt tr:knn:300} is the $k$NN predictor for transmission constraints. This predictor does not add any hyperplanes to the relaxation, and is used only as a baseline. All the following predictors were built on top of {\tt tr:knn:300}.
Predictor {\tt aff:svm} is the affine subspace predictor described in Subsection~\ref{subsec:aff}, which uses Support Vector Machines and cross-validation to find hyperplanes. Predictors {\tt aff:A} and {\tt aff:B} follow the same strategy, but use different hyperparameters. For {\tt aff:A}, we set $z^{\text{FIX}}=z^{\text{MAX}} = 0.975, z^{\text{MIN}}=0.025$ and keep all the other values unchanged. Our goal here is to allow prediction even under significant class imbalance; and when classes are too imbalanced, we simply predict the most common class. For {\tt aff:B} we set $z^{\text{FIX}}, z^{\text{MAX}}$ and $z^{\text{MIN}}$ to the same values as {\tt aff:B}, and, in addition, we lower all precision and recall thresholds $\alpha^R$ and $\alpha^P$ to 0.50, in an attempt to make {\tt aff:A} even more aggressive. Predictor {\tt aff:C} uses the same hyperparameters as {\tt aff:svm}, but we disable cross-validation. That is, if the classes are balanced enough, we simply train an SVM classifier and blindly use it to make predictions, without evaluating its accuracy. Finally, {\tt aff:perf} represents the perfect predictor, which fixes all binary variables to their optimal values, with perfect accuracy and zero running time. To implement this predictor, we gave it access to the optimal solution produced by {\tt tr:knn:300}.

Table~\ref{table:affCompare} shows the average wallclock running time (in seconds) for each predictor and the speedup relative to {\tt zero}. The table also shows, in column ``feasible'', how frequently the problem variations became infeasible after adding the hyperplanes. The three ``gap'' columns show the distribution of relative optimality gaps. For example, a value of 0.07 in the 95\% column means that, for 95\% of the variations, the method produced solutions with gap 0.07\% or better. The 100\% gap column shows the worst optimality gap observed for any variation.
\begin{table}
  \TABLE
  {Performance comparison of different affine subspace predictors. \label{table:affCompare}}
  {
    \setlength{\tabcolsep}{0.5em}
    \begin{tabular}{lrrrrrr}
        \toprule
        &
        &
        \multicolumn{3}{c}{Gap (\%)} &
        &
        \\
        \cmidrule(r){3-5}
        Method &
        Feasible (\%) &
        80\% &
        95\% &
        100\% &
        Time (s) &
        Speedup
        \\
        \midrule
        \begin{filecontents*}{tables/aff_compare.csv}
algorithm,wallclockTime,speedupZero,gap 0.800,gap 0.950,gap 1.000,success
tr:knn:300,80.36,2.82,0.04,0.07,0.09,100.00
aff:svm,22.24,10.19,0.03,0.05,0.08,100.00
aff:A,19.98,11.35,0.04,0.05,0.17,100.00
aff:B,12.04,18.83,0.05,0.11,0.25,72.59
aff:C,3.34,67.94,0.15,0.32,0.63,24.07
aff:perf,3.71,61.10,0.05,0.08,0.09,100.00
\end{filecontents*}
\csvreader[
                late after line={\\},
                late after last line={\\},
                filter expr={
                    test{\ifnumless{\thecsvinputline}{11}}
                }
            ]{tables/aff_compare.csv}{
                1=\algorithm,
                2=\wallclockTime,
                3=\speedupZero,
                4=\gapEighty,
                5=\gapNinetyFive,
                6=\gapHundred,
                7=\success,
            }{%
                \texttt{\algorithm} &
                \success &
                \gapEighty &
                \gapNinetyFive &
                \gapHundred &
                \wallclockTime &
                \speedupZero
            }%
        \bottomrule
    \end{tabular}
  }
  {}
\end{table}

Similarly to the prediction of warm starts, we observe that the prediction of valid hyperplanes is a significantly hard machine learning task. Blindly applying SVM to the problem, without filtering for low-quality predictors, had disastrous consequences. On average, {\tt aff:C} made 76\% of the test variations infeasible. Even for variations that remained feasible, the additional hyperplanes caused a significant deterioration in solution quality. Predictors {\tt aff:svm}, {\tt aff:A} and {\tt aff:B}, which use cross-validation to filter out low-quality SVM predictors had significantly better performance. Predictor {\tt aff:svm}, our recommended predictor, and the most conservative, presented an average speedup of 10.2x and never produced invalid or sub-optimal solutions during our tests. Predictor {\tt aff:A} also performed very well in the typical case, with an average speedup of 11.4x, however produced sub-optimal solutions for a small number of instance variations (less than 5\%), with the worst optimality gap reaching 0.17\%. Predictor {\tt aff:B} was overly aggressive, and, although not as often as {\tt aff:C}, still produced large number of infeasible problems and sub-optimal solutions. Finally, we see that, despite the impressive performance of {\tt aff:svm} and {\tt aff:A}, these predictors are still very far from the optimal speedup of 61x obtained by {\tt aff:perf}. We were unfortunately not able to reduce this gap any further.

In the following, we have a closer look the performance of {\tt aff:svm}. For each instance, Table~\ref{table:affine} shows the total number of commitment variables in the formulation, the percentage of commitment variables $x_{gt}$ fixed to zero, fixed to one, fixed to the next time period (that is, $x_{gt} = x_{g,t+1}$) by {\tt aff:svm}, as well as the percentage of commitment variables left free.
As explained in Subsection~\ref{subsec:aff}, each variable is fixed to at most one value, and therefore the percentages sum to 100\%.
\begin{table}
  \TABLE
  {Total number of commitment variables; percentage of fixed variables by affine predictor and percentage of correctly-identified hyperplanes. \label{table:affine}}
  {
    \setlength{\tabcolsep}{0.5em}
    \begin{tabular}{lrrrrrr}
        \toprule
        & \multicolumn{5}{c}{Commitment Variables}
        \\
        \cmidrule(r){2-6}
        Instance
        & Total
        & Fix Zero
        & Fix One
        & Fix Next
        & Free
        & Precision
        \\
        \midrule
        \begin{filecontents*}{tables/affine.csv}
instance,total,fixOne,fixZero,fixNext,free,precision
case1888rte,7128.0,52.0,17.9,24.7,5.4,99.8
case1951rte,9384.0,45.9,20.0,27.8,6.3,99.8
case2848rte,13128.0,44.8,23.9,23.2,8.2,99.7
case3012wp,12048.0,54.5,18.5,21.3,5.6,99.8
case3375wp,14304.0,57.8,14.4,22.4,5.5,99.8
case6468rte,31080.0,72.4,7.6,12.8,7.2,99.9
case6470rte,31920.0,70.3,7.3,16.2,6.2,99.8
case6495rte,32928.0,75.4,5.8,13.1,5.7,99.9
case6515rte,33312.0,75.5,6.1,12.0,6.4,99.9
Average,20581.3,60.9,13.5,19.3,6.3,99.8
\end{filecontents*}
\csvreader[
                head to column names,
                late after line={\\},
                late after last line={\\},
                filter expr={
                    test{\ifnumless{\thecsvinputline}{11}}
                }
            ]{tables/affine.csv}{}{%
                \texttt{\instance} &
                \total &
                \fixZero \% &
                \fixOne \% &
                \fixNext \% &
                \free \% &
                \precision \%
            }%
          \midrule
        \csvreader[
                head to column names,
                late after line={\\},
                late after last line={\\},
                filter expr={
                    test{\ifnumgreater{\thecsvinputline}{10}}
                }
            ]{tables/affine.csv}{}{%
                \textbf{\instance} &
                \total &
                \fixZero \% &
                \fixOne \% &
                \fixNext \% &
                \free \% &
                \precision \%
            }%
        \bottomrule
    \end{tabular}
  }
  {}
\end{table}
We observe that {\tt aff:svm} was able to fix a large number of commitment variables, for instances of all sizes. On average, the predictor was able to eliminate 94\% of the commitment variables, leaving only 6\% free, a considerable reduction in problem size. To evaluate the accuracy of such variable fixing decisions, we compared the recommendations provided by the predictor against an optimal solution obtained by the MIP solver (without machine learning). We note that such a comparison is not entirely fair to the predictor, since a large number of optimal solutions exist within the 0.1\% optimality gap threshold. These solutions have slight variations, and disagree with each other in a number of commitment variables.
Nevertheless, on average, 99.8\% of the hints provided by the predictor agreed the optimal solution obtained by CPLEX.

\begin{figure}
  \FIGURE
  {\includegraphics[width=0.75\textwidth]{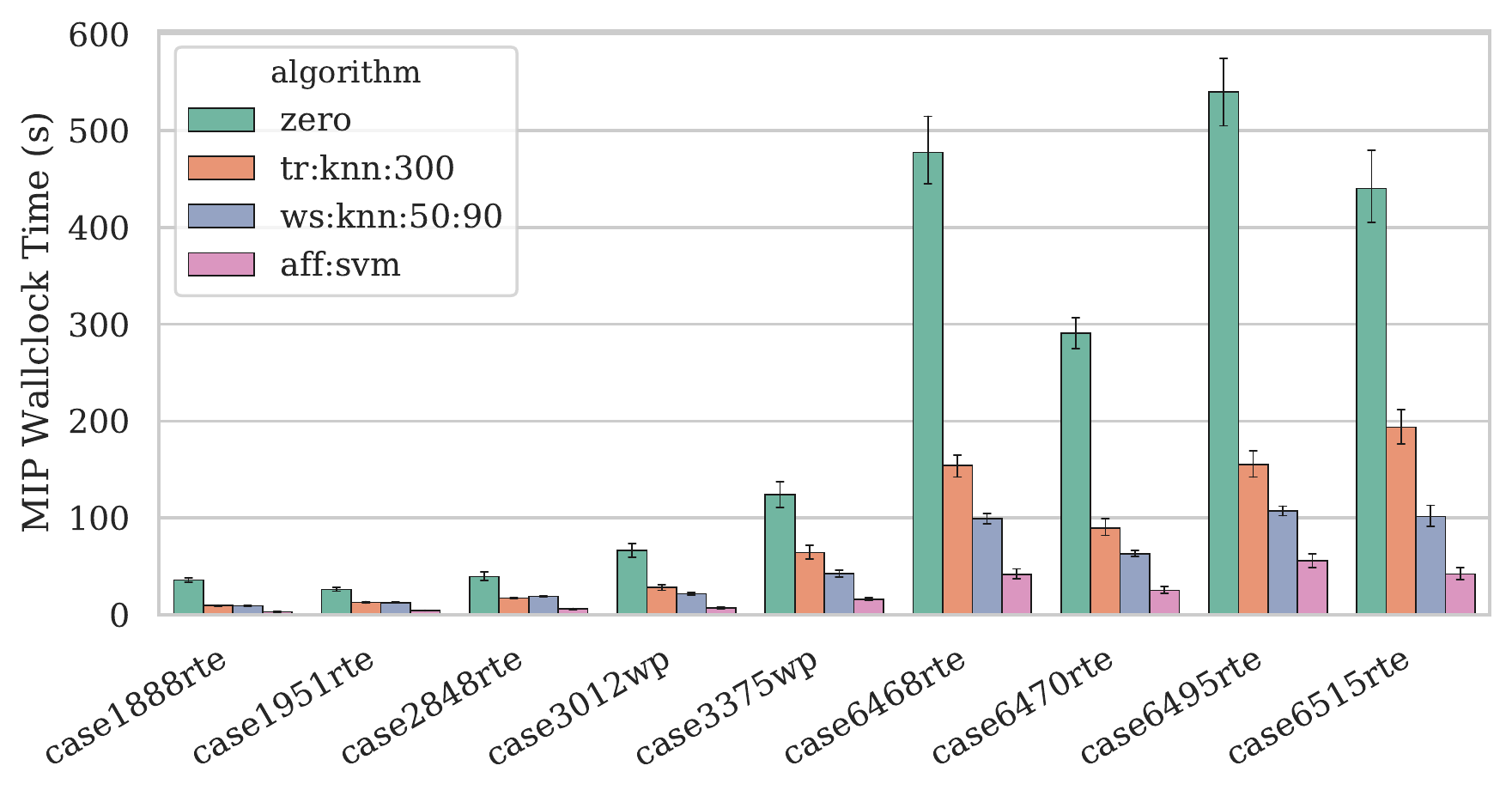}}
  {Running time. \label{fig:speed}}
  {}
\end{figure}

Table~\ref{table:speed} shows a summary of the results obtained.
Four different strategies are presented: no machine learning (\texttt{zero}), transmission prediction (\texttt{tr:knn:300}), transmission plus warm-start predicton (\texttt{ws:knn:50:90}) and finally transmission plus affine subspace prediction (\texttt{aff:svm}). Figure~\ref{fig:speed} shows the same results in graphical form, along with error bars representing 95\% confidence intervals. Overrall, by using different machine learning strategies, we were able to obtain a 4.3x speedup while maintaining optimality guarantees, and 10.2x speedup without guaratees, but with no observed reduction in solution quality.
\begin{table}
  \TABLE
  {Impact of machine-learning predictors on running-time (in-distribution). \label{table:speed}}
  {
    \setlength{\tabcolsep}{0.5em}
    \begin{tabular}{lrrrrrrr}
        \toprule
        & \multicolumn{1}{c}{\texttt{zero}}
        & \multicolumn{2}{c}{\texttt{tr:knn:300}}
        & \multicolumn{2}{c}{\texttt{ws:knn:50:90}}
        & \multicolumn{2}{c}{\texttt{aff:svm}}
        \\
        \cmidrule(r){2-2}
        \cmidrule(r){3-4}
        \cmidrule(r){5-6}
        \cmidrule(r){7-8}
        Instance
        &\footnotesize Time
        &\footnotesize Time
        &\footnotesize Speedup
        &\footnotesize Time
        &\footnotesize Speedup
        &\footnotesize Time
        &\footnotesize Speedup
        \\
        \midrule
        \begin{filecontents*}{tables/aff_speed.csv}
instance,wallclockTime:aff:svm,wallclockTime:tr:knn:300,wallclockTime:ws:knn:50:90,wallclockTime:zero,speedupZero:aff:svm,speedupZero:tr:knn:300,speedupZero:ws:knn:50:90,speedupZero:zero
case1888rte,2.83,9.26,9.15,35.68,12.61,3.85,3.90,1.00
case1951rte,4.19,12.68,12.50,26.07,6.22,2.06,2.09,1.00
case2848rte,5.96,17.16,18.75,39.38,6.61,2.29,2.10,1.00
case3012wp,7.03,27.96,21.48,66.11,9.41,2.36,3.08,1.00
case3375wp,15.88,64.17,42.24,124.05,7.81,1.93,2.94,1.00
case6468rte,41.76,154.04,99.42,477.48,11.43,3.10,4.80,1.00
case6470rte,25.14,89.54,63.10,290.55,11.56,3.24,4.60,1.00
case6495rte,55.52,154.92,107.12,540.38,9.73,3.49,5.04,1.00
case6515rte,41.88,193.56,101.47,440.37,10.52,2.28,4.34,1.00
Average,22.24,80.36,52.80,226.67,10.19,2.82,4.29,1.00
\end{filecontents*}
\csvreader[
                late after line={\\},
                late after last line={\\},
                filter expr={
                    test{\ifnumless{\thecsvinputline}{11}}
                }
            ]{tables/aff_speed.csv}{%
                1=\instance,
                2=\WAff,
                3=\WTr,
                4=\WWs,
                5=\WZero,
                6=\SAff,
                7=\STr,
                8=\SWs,
                9=\SZero,
            }{%
                \texttt{\instance} &
                \WZero&
                \WTr&
                \STr x &
                \WWs&
                \SWs x &
                \WAff &
                \SAff x
            }%
          \midrule
          \csvreader[
                  late after line={\\},
                  late after last line={\\},
                  filter expr={
                      test{\ifnumgreater{\thecsvinputline}{10}}
                  }
              ]{tables/aff_speed.csv}{%
                1=\instance,
                2=\WAff,
                3=\WTr,
                4=\WWs,
                5=\WZero,
                6=\SAff,
                7=\STr,
                8=\SWs,
                9=\SZero,
            }{%
                \textbf{\instance} &
                \WZero&
                \WTr&
                \STr x &
                \WWs&
                \SWs x &
                \WAff &
                \SAff x
            }%
        \bottomrule
    \end{tabular}
  }
  {}
\end{table}

\subsection{Out-of-Distribution Evaluation}
\label{subsec:ood}

In previous subsections, we evaluated the performance of the Machine Learning predictors under the assumption that test and training samples are drawn from exactly the same probability distribution.
This is a common assumption in the literature of predictive modeling, but it may not be entirely realistic, especially when there is a limited amount of historical data, or when the distribution shifts over time.
In this subsection, we evaluate the computational performance of the predictors when the test and training distributions are somewhat similar, but not exactly the same.
This scenario is known in the literature as \emph{dataset shift} \citep{Quionero2009}.
Our goal here is simply to evaluate whether the performance of the predictors dramatically deteriorates under moderate dataset shift, not to establish that the predictors are completely robust against such changes. We still expect the predictors to perform the best and most reliably when the training and test distributions are as similar as possible.

In the following, we use the same training set as before, where the 300 training samples are drawn from the distribution described in Subsection~\ref{subsec:samples}, but we now draw the test samples from a modified distribution, constructed as follows:
\begin{enumerate}[(i)]
    \item In the original distribution, production and startup cost multipliers are drawn from the uniform distribution in the interval $[0.95, 1.05]$. We replace this by the Gaussian distribution $\mathcal{N}(\mu, \sigma^2)$ with mean $\mu=1.05$ and standard deviation $\sigma=0.017$, so that multipliers now fall, with 99.7\% likelihood, in the interval $[1.00, 1.10]$, but may occasionally fall outside.
    That is, we assume that, during training, we correctly identified the size of the interval where multipliers are likely to fall, but incorrectly identified the mean, the distribution, and we did not consider the possibility of outliers.
    
    \item For geographical load distribution, we similarly replace the uniform distribution in the interval $[0.90, 1.10]$ by the Gaussian distribution $\mathcal{N}(\mu, \sigma^2)$ with mean $\mu=1$ and standard deviation $\sigma=0.033$, so that multipliers fall, with high likelihood, in the interval $[0.90, 1.10]$. Modifying the mean has no impact in this case, since load percentages are later normalized, as described in Subsection~\ref{subsec:samples}.
    Although the interval is the same, this distribution is still significantly different, and the test samples may now include outliers.
    
    \item For peak demand, we replace the uniform distribution in the interval $[0.555, 0.645]$ by the Gaussian distribution $\mathcal{N}(\mu, \sigma^2)$ with mean $\mu=0.60 \times 1.03$ and standard deviation $\sigma=0.009$, so that multipliers are likely to fall in the interval $[0.591, 0.645]$. That is, we assume that, during test, the mean peak demand is 3\% higher than during training and that the distribution is significantly different. We still assume that the multipliers fall, with very high likelihood, within an interval that was explored during training, but may now contain outliers outside this interval.
\end{enumerate}
The modifications above can be categorized as \emph{covariate shift}, since we only modify the probability of a certain instance variation being selected.
As explained in Section~\ref{sec:learning}, because we focus on the scenario where the training takes places the day before market clearing, we assume that training and test samples are solved using the same mathematical model (therefore, we have no \emph{prior probability shift}) and that new instances can still be encoded using the same set of parameters (therefore, we have no \emph{concept shift}). 

Table~\ref{table:speed-ood} 
summarizes the performance of predictors {\tt tr:knn:300}, {\tt ws:knn:50:90} and {\tt aff:svm} for this modified dataset. The table shows the average wallclock time (in seconds) necessary to solve each group of instance and the speedup relative to {\tt zero}.
Compared to Table~\ref{table:speed}, we observe that this moderate dataset shift caused some performance deterioration in all predictors, as expected, but nothing very extreme. The speedups of {\tt tr:knn:300}, {\tt ws:knn:50:90} and {\tt aff:svm} dropped from 2.8x, 4.3x and 10.2x, respectively, to 2.5x, 3.6x and 8.6x. Even with a degraded training set, the predictors still provided very strong performance benefits.

\begin{table}
  \TABLE
  {Impact of machine-learning predictors on running-time (out-of-distribution). \label{table:speed-ood}}
  {
    \setlength{\tabcolsep}{0.5em}
    \begin{tabular}{lrrrrrrr}
        \toprule
        & \multicolumn{1}{c}{\texttt{zero}}
        & \multicolumn{2}{c}{\texttt{tr:knn:300}}
        & \multicolumn{2}{c}{\texttt{ws:knn:50:90}}
        & \multicolumn{2}{c}{\texttt{aff:svm}}
        \\
        \cmidrule(r){2-2}
        \cmidrule(r){3-4}
        \cmidrule(r){5-6}
        \cmidrule(r){7-8}
        Instance
        &\footnotesize Time
        &\footnotesize Time
        &\footnotesize Speedup
        &\footnotesize Time
        &\footnotesize Speedup
        &\footnotesize Time
        &\footnotesize Speedup
        \\
        \midrule
        \begin{filecontents*}{tables/aff_speed_ood.csv}
instance,wallclockTime:aff:svm,wallclockTime:tr:knn:300,wallclockTime:ws:knn:50:90,wallclockTime:zero,speedupZero:aff:svm,speedupZero:tr:knn:300,speedupZero:ws:knn:50:90,speedupZero:zero
case1888rte,2.56,8.98,8.70,44.79,17.47,4.99,5.15,1.00
case1951rte,4.81,14.20,14.12,29.76,6.19,2.10,2.11,1.00
case2848rte,6.42,18.19,17.87,43.44,6.77,2.39,2.43,1.00
case3012wp,8.17,32.07,23.78,75.84,9.28,2.36,3.19,1.00
case3375wp,19.31,77.49,58.74,145.45,7.53,1.88,2.48,1.00
case6468rte,49.72,169.92,132.32,426.24,8.57,2.51,3.22,1.00
case6470rte,27.61,91.54,68.80,268.15,9.71,2.93,3.90,1.00
case6495rte,72.11,179.08,133.24,512.63,7.11,2.86,3.85,1.00
case6515rte,36.48,189.41,91.19,406.99,11.16,2.15,4.46,1.00
Average,25.24,86.76,60.97,217.03,8.60,2.50,3.56,1.00
\end{filecontents*}
\csvreader[
                late after line={\\},
                late after last line={\\},
                filter expr={
                    test{\ifnumless{\thecsvinputline}{11}}
                }
            ]{tables/aff_speed_ood.csv}{%
                1=\instance,
                2=\WAff,
                3=\WTr,
                4=\WWs,
                5=\WZero,
                6=\SAff,
                7=\STr,
                8=\SWs,
                9=\SZero,
            }{%
                \texttt{\instance} &
                \WZero&
                \WTr&
                \STr x &
                \WWs&
                \SWs x &
                \WAff &
                \SAff x
            }%
          \midrule
          \csvreader[
                  late after line={\\},
                  late after last line={\\},
                  filter expr={
                      test{\ifnumgreater{\thecsvinputline}{10}}
                  }
              ]{tables/aff_speed_ood.csv}{%
                1=\instance,
                2=\WAff,
                3=\WTr,
                4=\WWs,
                5=\WZero,
                6=\SAff,
                7=\STr,
                8=\SWs,
                9=\SZero,
            }{%
                \textbf{\instance} &
                \WZero&
                \WTr&
                \STr x &
                \WWs&
                \SWs x &
                \WAff &
                \SAff x
            }%
        \bottomrule
    \end{tabular}
  }
  {}
\end{table}

Regardless of the quality of the training dataset, we recall that {\tt tr:knn:300} and {\tt ws:knn:50:90} maintain all MIP optimality guarantees, and therefore can never negatively affect solution quality. The only potential negative impact of dataset shift, for these predictors, is a reduced speedup. For predictor {\tt aff:svm}, on the other hand, a sufficiently large dataset shift can potentially compromise both speed and solution quality.
To evaluate the impact of dataset shift in the quality of the solutions produced by {\tt aff:svm}, we present in Table~\ref{table:gap-ood} the distribution of their relative optimality gaps. A value of 0.08 in column 95\%, for example, means that, for 95\% of the variations of that instance, the relative optimality gap was {\tt aff:svm} was 0.08\% or better. Column 100\% shows the worst optimality gap obtained for any variation of that instance.

As Table~\ref{table:gap-ood} shows, even under the moderate dataset shift considered, {\tt aff:svm} still produced optimal solutions in the vast majority of the cases. For 8 out of the 9 instances considered, and particularly for all instances with more than 6000 buses, {\tt aff:svm} produced optimal solutions for all variations. For instance {\tt case3012wp}, around 5\% of the variations were sub-optimal, with the worst relative optimality gap reaching 0.13\%, which is only slightly higher than the 0.1\% threshold.

\begin{table}
  \TABLE
  {Solution quality (out-of-distribution). \label{table:gap-ood}}
  {
    \setlength{\tabcolsep}{0.5em}
    \begin{tabular}{lrrrr}
        \toprule
        & \multicolumn{4}{c}{Gap (\%)}
        \\
        \cmidrule(r){2-5}
        Instance
        &\footnotesize 50\%
        &\footnotesize 80\%
        &\footnotesize 95\%
        &\footnotesize 100\%
        \\
        \midrule
        \begin{filecontents*}{tables/gap_ood.csv}
instance,0.50,0.80,0.95,1.00
case1888rte,0.03,0.04,0.08,0.09
case1951rte,0.02,0.03,0.04,0.05
case2848rte,0.03,0.04,0.08,0.10
case3012wp,0.06,0.09,0.10,0.13
case3375wp,0.03,0.04,0.05,0.06
case6468rte,0.02,0.03,0.07,0.07
case6470rte,0.03,0.04,0.05,0.07
case6495rte,0.02,0.03,0.05,0.05
case6515rte,0.02,0.03,0.04,0.05
\end{filecontents*}
\csvreader[
                late after line={\\},
                late after last line={\\},
                filter expr={
                    test{\ifnumless{\thecsvinputline}{11}}
                }
            ]{tables/gap_ood.csv}{%
                1=\instance,
                2=\fifty,
                3=\eighty,
                4=\ninetyfive,
                5=\hundred
            }{%
                \texttt{\instance} &
                \fifty &
                \eighty &
                \ninetyfive &
                \hundred
            }%
        \bottomrule
    \end{tabular}
  }
  {}
\end{table}

\section{Limitations and Future Work} \label{sec:limitations}

In this work, we proposed three Machine Learning predictors for expediting the solution of the Security-Constrained Unit Commitment Problem, and evaluated their performance on a number of realistic, large-scale instances of the problem. By predicting redundant constraints and warm starts, we obtained a 4.3x speedup over the baseline, while still maintaining all optimality guarantees. By predicting a restricted affine subspace where the solution was very likely to reside, we increased this speedup to 10.2x. Although no optimality guarantees were produced in this case, we observed that all solutions produced in our experiments were, in fact, optimal. We also performed out-of-distribution experiments, where the training and test distribution did not exactly match. Even under this challenging scenario, the predictors presented very strong performance, indicating that they are at least somewhat robust against dataset shift.

The main limitation of our approach is that, in order to train the predictors, a large number of solved instances must be available, and they need to be sufficiently similar to the instances we expect to solve in the future.
In our experiments, we considered the scenario where training instances are artificially generated and solved the day before market clearing takes places. Sufficient historical data must be available to characterize the distribution of parameters that are uncertain the day before market clearing, so that plausible training instances can be generated. There is also some upfront computational cost when solving the training instances, although this limitation is mitigated by the fact that: (i) training instances can be solved in parallel, and the number of samples required by the proposed method is relatively small; (ii) the first two predictors presented in this work can be used in online fashion, and therefore can be used to significantly accelerate the generation of the training dataset itself.
Future work will consider ML methods that are able to handle changes to the generation fleet, so that training can be performed less frequently.
Alternatively, we are exploring methods to update an existing training dataset more efficiently, instead of recreating it from scratch, following a significant dataset shift.
Further validation on real-world datasets is also needed.
Finally, we note that the techniques presented here can be adapted to other challenging combinatorial problems.



\newcommand{\ACKNOWLEDGMENT}{\section*{Acknowledgments}}
\newcommand{\APPENDIX}{\section*{Appendix}}

\ACKNOWLEDGMENT{%
This material is based upon work supported by Laboratory Directed Research and Development (LDRD) funding from Argonne National Laboratory, provided by the Director, Office of Science, of the U.S. Department of Energy under Contract No. DE-AC02-06CH11357.
We gratefully acknowledge use of the Bebop cluster in the Laboratory Computing Resource Center at Argonne National Laboratory.
}%

\APPENDIX{
Here we present the complete MIP formulation that was used to during the computational experiments in Section~\ref{sec:comp}.
Consider a power system composed by a set $B$ of buses, a set $G$ of generators and a set $L$ of transmission lines. Furthermore, let $T=\{1,\ldots,24\}$ be the set of hours within the planning horizon, and let $G_b$ be the set of generators located at bus $b$.
Let $d_{bt}$ be the demand (in MWh) from bus $b$ at time $t$. We recall that each generator $g \in G$ has a convex production cost curve, modeled as a piecewise-linear function with a set $K$ of segments. For each generator $g \in G$, define the following constants:
\begin{itemize}
  \item $c^0_g$, cost of keeping generator operational for one hour, producing at its minimum output level.
  \item $c^k_g$, cost to produce each additional MWh of power, for each segment $k \in K$.
  \item $c^S_g$, cost to start the generator up.
  \item $P^k_g$, amount of power available (in MWh) in segment $k \in K$.
  \item $RU_g$, maximum allowed rise in production (in MWh) from one hour to the next.
  \item $RD_g$, maximum allowed drop in production (in MWh) from one hour to the next.
  \item $P^{\min}_g$, minimum amount of power (in MWh) the generator must produce if it is operational.
  \item $P^{\max}_g$, maximum amount of power (in MWh) the generator can produce.
  \item $UT_g$, minimum amount of time (in hours) the generator must stay operational after being switched on.
  \item $DT_g$, minimum amount of time (in hours) the generator must stay offline after being switched off.
\end{itemize}
For each transmission line $l \in L$, let $F^0_l$ be its flow limit (in MWh) under normal conditions and $F^c_l$ be its flow limit when there is an outage on line $c \in L$. Similarly, let $\delta^0_{lb}$ and $\delta^c_{lb}$ be, respectibely, the injection shift factors for line $l$ and bus $b$ under normal conditions, and under outage on line $c$.

As mentioned in Subsection~\ref{subsec:scuc}, the main decision variables for this problem are $x_{gt} \in \{0,1\}$, which indicates whether generator $g$ is operational at time $t$; and $y_{gt} \geq 0$, which indicates how much power (in MWh) generator $g$ produces during time $t$.
Other auxiliary variables are $z_{gt}, w_{gt} \in \{0,1\}$ which indicate, respectively, whether generator $g$ is started up or shut down at time $t$.
We also define the variables $y^k_{gt}$ which indicate how much power produced by generator $g$ at time $t$ comes from segment $k \in K$.
Finally, let $r_{gt} \geq 0$ be a decision variable indicating the amount of reserve (in MWh) provided by generator $g$ at time $t$, and let $R_t$ be the minimum amount of system-wide reserve required at time $t$. Reserves are generation capacities kept aside to compensate for small load variations.
Given these variables and constants, SCUC is formulated as a cost minimization problem as below:

\small
\begin{align}
\text{Minimize} \hspace{2em} &
    \sum_{g \in G} \sum_{t \in T} \left[
        c^S_g z_{gt} + c^0_g x_{gt} + \sum_{k \in K} c^k_g y^k_{gt}
    \right]
    \label{eq:obj} \\
\text{Subject to} \hspace{2em} &
    \sum_{g \in G} y_{gt} = \sum_{b \in B} d_{bt}
        & \forall t \in T \label{eq:balance}\\
    & \sum_{g \in G} r_{gt} \geq R_t
        & \forall t \in T \label{eq:reserve}\\
    & \sum_{k \in K} y^k_{gt} + r_{gt} \leq
        (P^{\max}_g - P^{\min}_g) x_{gt} \notag\\
        & \hspace{6em} - (P^{\max}_g - RU_g) z_{gt} \notag\\
        & \hspace{6em} - (P^{\max}_g - RD_g) w_{g,t+1}
        & \forall g \in G, t \in \{1,\ldots,23\} : UT_g > 1
        \label{eq:capacity1}\\
    & \sum_{k \in K} y^k_{gt} + r_{gt} \leq
        (P^{\max}_g - P^{\min}_g) x_{gt} \notag\\
        & \hspace{6em} - (P^{\max}_g - RU_g) z_{gt} \notag \\
        & \hspace{6em} - \max\{RU_g - RD_g, 0\} w_{g,t+1}
        & \forall g \in G, t \in \{1,\ldots,23\} : UT_g = 1
        \label{eq:capacity2} \\
    & \sum_{k \in K} y^k_{gt} + r_{gt} \leq
        (P^{\max}_g - P^{\min}_g) x_{gt} \notag\\
        & \hspace{6em} - (P^{\max}_g - RD_g) w_{g,t+1} \notag\\
        & \hspace{6em} - \max\{RD_g - RU_g, 0\} z_{gt}
        & \forall g \in G, t \in \{1,\ldots,23\} : UT_g = 1
        \label{eq:capacity3} \\
    & \sum_{k \in K} y^k_{g,24} + r_{g,24} \leq
        (P^{\max}_g - P^{\min}_g) x_{g,24} \notag\\
        & \hspace{8em} - (P^{\max}_g - RU_g) z_{t,24}
        & \forall g \in G
        \label{eq:capacity4} \\
    & y_{gt} \leq y_{g,t-1} + RU_g
        & \forall g \in G, t \in \{2,\ldots,24\} \label{eq:ru} \\
    & y_{gt} \geq y_{g,t-1} - RD_g
        & \forall g \in G, t \in \{2,\ldots,24\} \label{eq:rd} \\
    & \sum_{i=\max\{1,t-UT_g+1\}}^t z_{gi} \leq x_{gt}
        & \forall t \in T, g \in G \label{eq:ut} \\
    & \sum_{i=\max\{1,t-DT_g+1\}}^t z_{gi} \leq 1 - x_{\max\{1,t-DT_g\},g}
        & \forall t \in T, g \in G \label{eq:dt} \\
    & -F^0_l \leq \sum_{b \in B} \delta^0_{lb} \left( \sum_{g \in G_b} y_{gt} - d_{bt} \right) \leq F^0_l & \forall l \in L, t \in T \label{eq:tr1} \\
    & -F^c_l \leq \sum_{b \in B} \delta^c_{lb} \left( \sum_{g \in G_b} y_{gt} - d_{bt} \right) \leq F^c_l & \forall c \in L, l \in L, t \in T \label{eq:tr2} \\
    & y^k_{gt} \leq P^k_g & \forall k \in K, g \in G, t \in T \label{eq:seg1} \\
    & y_{gt} = P^{\min}_g x_{gt} + \sum_{k \in K} y^k_{gt}
        & \forall t \in T \label{eq:seg2} \\
    & x_{gt} - x_{g,t-1} = z_{gt} - w_{gt}
        & \forall g \in G, t \in \{2,\ldots,24\} \label{eq:link} \\
    & x_{gt}, z_{gt}, w_{gz} \in \{0,1\} & \forall g \in G, t \in T \\
    & r_{gt}, y_{gt}, y^k_{gt} \geq 0 & \forall k \in K, g \in G, t \in T
\end{align}
The objective function \eqref{eq:obj} includes start-up and production costs. Although start-up costs are sometimes modeled as a stepwise function of off-time, in our test they are modeled as constants.
Equation \eqref{eq:balance} enforces that the total power supply equals total load.
Equation \eqref{eq:reserve} enforces a sufficient amount of reserve at each time period.
Equations \eqref{eq:capacity1} to \eqref{eq:capacity4} enforce the production limits.
Equations \eqref{eq:ru} and \eqref{eq:rd} enforce the ramping requirements.
Equations \eqref{eq:ut} and \eqref{eq:dt} guarantee that, once a generator is started or shutdown, it must remain online or offline for a certain amount of time.
Equations \eqref{eq:tr1} and \eqref{eq:tr2} require that the power flow on each transmission line does not exceed its thermal limits.
Equations \eqref{eq:seg1} and \eqref{eq:seg2} link the variables $y_{gt}$ and $y^k_{gt}$.
Finally, Equation \eqref{eq:link} link $x_{gt}, z_{gt}$ and $w_{gt}$.
}

\bibliographystyle{apalike}%
\bibliography{papers.bib}%

\end{document}